\documentclass[doublespace,12pt]{amsart}
\usepackage{bibunits}
\usepackage[english]{babel}
\usepackage[T1]{fontenc}
\usepackage[utf8]{inputenc}
\usepackage{amsmath,amssymb,amsthm}
\usepackage{textcomp}
\usepackage{stmaryrd}
\usepackage[sc]{mathpazo}
\usepackage[left=2.5cm,right=2.5cm,top=2cm,bottom=3cm]{geometry}
\usepackage{eulervm}
\usepackage[cal=dutchcal,scr=kp]{mathalpha}
\makeatletter
\def\mathalfa@frakscaled{s*[1.0]}
  \DeclareFontFamily{U}{esstixfrak}{\skewchar \font =45}
  \DeclareFontShape{U}{esstixfrak}{m}{n}{
    <-> \mathalfa@frakscaled esstixfrak}{}
  \DeclareMathAlphabet{\mathfrak}{U}{esstixfrak}{m}{n}
\makeatother

\usepackage{mathtools}
\usepackage{comment}
\usepackage{todonotes}
\setuptodonotes{textcolor=blue, linecolor=blue, bordercolor=blue, backgroundcolor=white, size=footnotesize}

\usepackage{hyperref}
\usepackage{tikz}
\usetikzlibrary{knots,hobby,decorations.pathreplacing,shapes.geometric,calc,decorations.text,shapes.misc,calc,math,arrows,decorations.markings,patterns, matrix, intersections,cd}
\usepackage[nameinlink]{cleveref}

\definecolor{couleur3}{rgb}{0.0, 0.5, 1.0} 
\definecolor{couleur4}{rgb}{0.55, 0.71, 0.0} 
\definecolor{couleur1}{rgb}{0.87, 0.45, 1.0} 
\definecolor{couleur2}{rgb}{0.81, 0.06, 0.13}
\definecolor{couleur5}{rgb}{1, .64, 0.28}

\newcommand{\catname}[1]{\mathbf{#1}}
\newcommand{\bicatname}[1]{\mathfrak{#1}}
\newcommand{\sitename}[1]{\mathfrak{#1}}
\newcommand{\nom}[1]{\textsc{#1}}

\newcommand{\B}{\mathbb{B}}
\newcommand{\C}{\mathbb{C}}

\newcommand{\N}{\mathbb{N}}

\renewcommand{\P}{\mathbb{P}}
\newcommand{\Q}{\mathbb{Q}}
\newcommand{\R}{\mathbb{R}}
\renewcommand{\S}{\mathbb{S}}
\newcommand{\T}{\mathbb{T}}
\newcommand{\U}{\mathbb{U}}

\newcommand{\Z}{\mathbb{Z}}

\newcommand{\Ccal}{\mathcal{C}}

\newcommand{\Fcal}{\mathcal{F}}

\newcommand{\Ical}{\mathcal{I}}

\newcommand{\pcal}{\mathcal{p}}

\newcommand{\Ascr}{\mathscr{A}}

\newcommand{\Escr}{\mathscr{E}}
\newcommand{\Fscr}{\mathscr{F}}

\newcommand{\Kscr}{\mathscr{K}}

\newcommand{\Mscr}{\mathscr{M}}

\newcommand{\Sscr}{\mathscr{S}}
\newcommand{\Tscr}{\mathscr{T}}
\newcommand{\Uscr}{\mathscr{U}}

\newcommand{\Xscr}{\mathscr{X}}

\newcommand{\Sfrak}{\mathfrak{S}}

\renewcommand{\hbar}{\overline{h}}

\newcommand{\hcal}{\mathcal{h}}
\newcommand{\lcal}{\mathcal{l}}

\DeclareMathOperator{\Hom}{\mathrm{Hom}}

\DeclareMathOperator{\Vect}{\mathrm{Vect}}

\DeclareMathOperator{\GL}{\mathrm{GL}}

\DeclareMathOperator{\im}{\mathrm{im}}
\DeclareMathOperator{\Cone}{\mathrm{Cone}}

\DeclareMathOperator{\Pic}{\mathrm{Pic}}

\DeclareMathOperator{\Aut}{\mathrm{Aut}}
\DeclareMathOperator{\comb}{\mathrm{comb}}
\DeclareMathOperator{\Gr}{\mathrm{Gr}}

\DeclareMathOperator{\Conv}{\mathrm{Conv}}

\crefname{subsection}{the subsection}{the subsections}
\crefname{section}{section}{the sections}
\crefname{subsubsection}{subsubsection}{the subsubsections}

\crefname{Thm}{theorem}{theorems}
\Crefname{Thm}{Theorem}{Theorems}
\crefname{Prop}{proposition}{propositions}
\Crefname{Prop}{Proposition}{Propositions}
\crefname{Lemme}{lemma}{lemmas}
\Crefname{Lemme}{Lemma}{Lemmas}
\crefname{Cor}{corollary}{corollaries}
\Crefname{Cor}{Corollary}{Corollaries}
\crefname{Const}{construction}{constructions}
\Crefname{Const}{Construction}{Constructions}
\crefname{Ex}{exemple}{exemples}
\Crefname{Ex}{Exemple}{Exemples}
\crefname{Def}{definition}{definitions}
\Crefname{Def}{Definition}{Definitions}
\crefname{Not}{notation}{notations}
\Crefname{Not}{Notation}{Notations}
\crefname{Rem}{remark}{remarks}
\Crefname{Rem}{Remarque}{Remarques}
\usepackage{cancel}

\begin{document}
\title{Gluing moduli spaces of quantum toric stacks via secondary fan}
\author{Antoine BOIVIN}
\email{\href{mailto:antoine.boivin@univ-angers.fr}{antoine.boivin@univ-angers.fr}}


\begin{abstract}
The extension from toric varieties to quantum toric stacks allows for the study of moduli spaces of toric objects with fixed combinatorial structures, as we now consider general finitely generated subgroups of $\R^n$ as "lattices." This paper aims to construct a moduli space that encompasses all such moduli spaces for a given dimension of the ambient space. To achieve this, we adapt the construction of the secondary fan within the quantum framework. This approach provides a description of wall-crossings between different moduli spaces, analogous to those observed in LVMB manifolds.
\end{abstract}

\date\today
\maketitle

\newtheorem{Thm}{Theorem}[subsection]
\newtheorem*{unThm}{Theorem}

\newtheorem{Prop}[Thm]{Proposition}

\newtheorem{Lemme}[Thm]{Lemma}

\newtheorem{Cor}[Thm]{Corollary}

\newtheorem{Const}[Thm]{Construction}

\newtheorem{DefProp}[Thm]{Definition-Proposition}

\newtheorem*{Quest}{Question}

\theoremstyle{definition}

\newtheorem{Ex}[Thm]{Example}

\newtheorem{Def}[Thm]{Definition}

\newtheorem{Not}[Thm]{Notation}

\theoremstyle{remark}
\newtheorem{Rem}[Thm]{Remark}
\newtheorem{Avert}[Thm]{Warning}

\setcounter{secnumdepth}{4}
\setcounter{tocdepth}{3}

\section*{Introduction}

Quantum toric geometry (introduced in \cite{katzarkov:hal-01672716}, \cite{katzarkov2020quantum} and \cite{boivin2020nonsimplicial}) is a generalization of toric geometry (as described in \cite{cox} or \cite{fulton}) which permits one to obtain moduli spaces of toric objects. For this purpose, the authors of \cite{katzarkov2020quantum} introduce "quantum toric stacks", which are analytic stacks built from the combinatorial data of the "quantum fan", i.e. the data of a good family of cones generated by elements of an arbitrary finitely generated subgroup $\Gamma$ of $\R^n$. This description is functorial and induces an equivalence between the category of quantum fans and the category of quantum toric stacks. The transition from the discrete to this general case allows us to continuously deform our fans and hence obtain continuous families of quantum toric stacks.

These moduli spaces are studied in \cite[section 11]{katzarkov2020quantum} and \cite{boivin2023moduli}. In the latter paper, the considered moduli spaces are those of quantum toric stacks associated with quantum fans with fixed combinatorics and number of generators of the group $\Gamma$. They are (differentiable) orbifolds obtained as the quotient of an open connected semi-algebraic subset by a finite group determined by the combinatorics. Furthermore, \cite{boivin2023moduli} establishes that these moduli spaces admit a natural compactification by incorporating quantum fans with degenerated combinatorics. 

The goal of this paper is to pursue this work by constructing a big moduli space of quantum toric stacks by connecting the different moduli spaces of the last paragraph. The main tool for this is (an adaptation of) the secondary fan in the quantum framework (see \cite[section 5]{Sottile} for another use of the secondary fan in the irrational setting). In the classical case (as introduced by \cite{gelfand2009discriminants}), if we have a (reductive) subgroup $G$ of a torus $(\C^*)^n$ and $\chi$ a character of $G$, then we can form the GIT quotient $\C^d//_\chi G$ which is a toric variety. The secondary fan is a partition of a cone of $\widehat{G} \otimes \R$ of the different combinatorics of the fan of $\C^d//_\chi G$ obtained in this way. Unfortunately, in the quantum framework, the group of the quotient is $\C^{n-d}$ which is not reductive. We will adapt in \cref{secondary_fan} a combinatoric construction suitable to the quantum case. It permits obtaining, in the general case, results of wall-crossing between different moduli spaces (using the formalism of \cite{boivin2024birational} of combinatorial birational morphism) and an analog of cobordism between LVMB manifolds (as described in \cite{bosio2006}) for quantum toric stacks.

In \cref{augmented_ms}, we will describe the universal family over this moduli space (said "augmented") extending the one on each smaller moduli space thanks to the wall-crossings encoded by the secondary fan.     

Throughout this article, we proved a strong result of connectedness of this big moduli space: 

\begin{unThm}[\cref{Thm_cnx_proj}]
    Let $(\Delta,h,\Ical)$ be a simplicial quantum fan. There exists a continuous path in the augmented moduli space between this quantum fan and the fan of a quantum projective space.
\end{unThm}

In conclusion, we get the main result of this paper:

\begin{unThm}[\cref{Compactification_edm_recolle}]
    Let $n \geq d$ two integers. There exists a compact stack $\Kscr(d,n)$ and a stack morphism $\overline{\Escr(d,n)} \to \Kscr(d,n)$ such that :
    \begin{itemize}
        \item For every combinatorial type $D$ of a complete fan of $\R^d$ with $n$ generators,there exists a closed substack $\Kscr(D)$ of $\Kscr(d,n)$ and a $\R^{n-d}$-fibration $\Kscr(D) \to \overline{\Mscr(d,n,D)}$ ;
        \item This fibration extends to a fibration between the family over $\Kscr(D)$ and the family over $\overline{\Mscr(d,n,D)}$. 
    \end{itemize}
    where $\overline{\Mscr(d,n,D)}$ is the compactification of the moduli space of quantum toric stacks of dimension $d$, with $n$ generators and of combinatorial type $D$.
    \end{unThm}

\textbf{Acknowledgement :} The author was supported by the Région Pays de la Loire from the France 2030 program, Centre Henri Lebesgue ANR-11-LABX-0020-01.

\tableofcontents
\section{Recall on quantum toric geometry and notations}
\subsection{Quantum fans and quantum toric stacks}
In this subsection, we recall the needed definitions and theorems on quantum toric stacks (see \cite{katzarkov2020quantum} for the details of the constructions).

In the classical case, we consider cones generated by elements of a lattice i.e. a discrete subgroup of rank $n$ of $\R^n$. In the quantum case, we want to consider general finitely generated subgroup of $\R^n$. Such group can have an arbitrary large rank and this rank is not constant by small deformations (for instance, the group $\Gamma_\alpha=\Z + \alpha \Z \subset \R$ is of rank $1$ if $\alpha \in \Q$ and of rank $2$ otherwise). In order to study moduli spaces in quantum toric geometry, we will fix the numbers of generators of the group:

\begin{Def}
Let $\Gamma$ be a finitely generated subgroup of $\R^d$ such that $\Vect_\R(\Gamma)=\R^d$. A calibration of $\Gamma$ is given by:
\begin{itemize}
\item A group epimorphism $h \colon \Z^n \to \Gamma$
\item A subset $\Ical \subset \{1,\ldots,n\}$ such that $\Vect_\C(h(e_j), j \notin \Ical)=\C^d$ (this is the set of virtual generators)

\end{itemize}
This is a standard calibration if $\Z^d \subset \Gamma$, $h(e_i)=e_i$ for $i=1,\ldots,d$ and $\Ical$ is of the form $\{n-|\Ical|+1,\ldots,n\}$ (it is not a restrictive condition (see \cite[Lemma 3.13]{katzarkov2020quantum})).
\end{Def}

\begin{Def} \label{calib_q_fan_def}
A quantum fan $(\Delta,h : \Z^n \to \Gamma \subset \R^d,\Ical)$ in $\Gamma$ is the data of 
\begin{itemize}
\item a collection $\Delta$ of strongly convex polyhedral cones generated by elements of $\Gamma$ such that every intersection of cones of $\Delta$ is a cone of $\Delta$, every face of a cone of $\Delta$ is a cone and $\{0\}$ is a cone of $\Delta$.
\item a standard calibration $h$ with $\Ical$ its set of virtual generators
\item A set of generators $A$ i.e. a subset of $\{1,\ldots,n\} \setminus \Ical$ such that the 1-cone generated by the $h(e_i)$ for $i \in A$ are exactly the 1-cones of $\Delta$ 
\end{itemize}  
The fan is said simplicial if every cone of $\Delta$ is simplicial (i.e. which can be send on a cone $\Cone(e_1,\ldots,e_k)$ (where $k$ is the dimension of the considered cone) by a linear automorphism of $\R^d$). \\
We note $\Delta(1)$ the cones of dimension 1 of $\Delta$ and $\Delta_{max}$ the maximal (for the inclusion) cones of $\Delta$. \\
With pair of linear morphisms $(L,H)$ which preserves inclusion of cones and calibrations (i.e. $(L,H) : (\Delta,h,\Ical) \to (\Delta',h',\Ical)$ then $hL=Hh'$), the quantum fans form a category denoted $\mathbf{QF}$. 
We note $\mathbf{QF}^{simp}$ the full subcategory of $\mathbf{QF}$ with simplicial fans.
\end{Def}

In \cite{katzarkov2020quantum}, the authors give a construction of a quantum toric stack associated to a simplicial quantum fan. It is a stack over the equivariant analytic site $\sitename{A}$\footnote{It is the category of analytic spaces with an action of an abelian complex Lie group with equivariant holomorphic maps and equivariant Euclidean covering.} defined by gluing of affine pieces as follows:

 Let $\sigma$ be a cone of a simplicial quantum fan $(\Delta,h : \R^n \to \R^d,\Ical)$. Let $L$ be a linear isomorphism of $\R^d$ such that $L(\sigma)=\Cone(e_1,\ldots,e_{\dim(\sigma)})$ and note $H$ the linear morphism described by a permutation $\chi$ of $\{1,\ldots,n\}$ such that $\chi(i_k)=k$ for $1 \leq k \leq \dim(\sigma)$). Then the cone $\sigma$ describe an affine quantum toric stack
\[
\Uscr_\sigma \coloneqq [\C^{\dim(\sigma)}\times (\C^*)^{d-\dim(\sigma)}/\Z^{n-d}]
\]

where the action of $\Z^{n-d}$ on $\C^{\dim(\sigma)}\times (\C^*)^{d-\dim(\sigma)}$ is given by the following morphism 
\[
x \in \Z^{n-d} \mapsto  EL h H^{-1}(0_{\R^{d}} \oplus x) \in (\C^*)^n
\]
where $E$ is the map \begin{equation} \label{def_exp}
    (z_1,\ldots,z_n) \in \C^n \mapsto (\exp(2i\pi z_1),\ldots,\exp(2i\pi z_n)) \in (\C^*)^n
\end{equation}
In particular, with $\sigma=\{0\}$, we get the quantum torus $\Tscr_{h,\Ical}\coloneqq[(\C^*)^d/\Z^{n-d}]$ (which is a Picard stack for the induced group law on the quotient) which is dense in all these quantum toric stacks. 

A Picard stack morphism $\mathcal{l}\colon\Tscr_{h,\Ical} \to \Tscr_{h',\Ical'}$  is combinatorially described by the data of two linear morphisms $(L : \R^d \to  \R^{d'},H : \R^n \to \R^{n'})$ such that the following diagram
\begin{align}\label{torus_morphism}
\begin{tikzcd}[ampersand replacement=\&]
	{\R^n} \& {\R^{n'}} \\
	{\R^d} \& {\R^{d'}}
	\arrow["H", from=1-1, to=1-2]
	\arrow["h"', from=1-1, to=2-1]
	\arrow["{h'}", from=1-2, to=2-2]
	\arrow["L"', from=2-1, to=2-2]
\end{tikzcd}
\end{align}
and which respects virtual generators (see \cite[Proposition 3.3.5]{Boivin}).

Thanks to the stability by taking faces and stability by intersection, a quantum fan $(\Delta,h,\Ical)$ define a diagram $D : \sigma \in \Sigma \mapsto \Uscr_\sigma$ of affine quantum toric stacks (where $\Delta$ is ordered by inclusion). Indeed, for every $\sigma,\tau \in \Delta$, we have
\begin{align}
\begin{tikzcd}[ampersand replacement=\&,scale=0.05]
    \&\& {U_\sigma} \\
	{U_{\sigma\cap\tau}} \\
	\&\& {U_\tau}
	\arrow["{\text{open }}"{pos=0.5}, hook, from=2-1, to=1-3]
	\arrow["{\text{open}}"'{pos=0.5}, hook, from=2-1, to=3-3]
\end{tikzcd}
\end{align}

\begin{Def}
    The quantum toric stack associated to the simplicial fan $(\Delta,h,\Ical)$ is the colimit $\Xscr_{\Delta,h,\Ical}$ of this diagram i.e. the gluing of the quantum toric stacks $\{\Uscr_\sigma\}_{\sigma \in \Delta}$ along their intersection. With stack morphisms which restrict on Picard stack morphism on tori, they form a category denoted $\mathbf{QTS}^{\mathrm{simp}}$.
\end{Def}

\begin{Thm} \label{equiv_cat}
    The correspondence $(\Delta,h,\Ical) \in \mathbf{QF}^{\mathrm{simp}} \mapsto \Xscr_{\Delta,h,\Ical} \in \mathbf{QTS}^{\mathrm{simp}}$ is an equivalence of categories.
\end{Thm}
This theorem permits us to study the moduli spaces of quantum toric stacks (see the next subsection)  with through the moduli spaces of quantum fans.

\begin{Thm}[Quantum GIT] \label{QGIT}
If $(\Delta,h,\Ical)$ is a simplicial quantum fan then the quantum toric stack $\Xscr_{\Delta,h,\Ical}$ is a quotient stack
\[
\Xscr_{\Delta,h,\Ical}=[\Sscr(\Delta)/\C^{n-d}]
\]
where $\Sscr(\Delta)$ is a quasi-affine (classical) toric variety given by the combinatorics of $\Delta$ :
    \begin{equation} \label{Def_S}
        \Sscr(\Delta)=\bigcup_{\sigma \in \Delta_{max}} \C^{\sigma(1)} \times (\C^*)^{\sigma(1)^c} \subset \C^n;
    \end{equation}    
and $\C^{n-d}$ acts on $\Sscr$ through 
\[
t \cdot z=E(k \otimes id_\C(t))z
\]
where $k$ is a Gale transform of $h \otimes id_\R : \R^n \to \R^d$ i.e. a linear map such that the short sequence 
\begin{align}\label{Gale_transform}
\begin{tikzcd}[ampersand replacement=\&]
	0 \& {\R^{n-d}} \& {\R^n} \& {\R^d} \& 0
	\arrow[from=1-1, to=1-2]
	\arrow["k", from=1-2, to=1-3]
	\arrow["h", from=1-3, to=1-4]
	\arrow[from=1-4, to=1-5]
\end{tikzcd}
\end{align}
is exact
\end{Thm}

\subsection{Moduli spaces of quantum toric stacks}

Thanks to the \cref{equiv_cat}, we can define moduli spaces of quantum toric stacks (see \cite{boivin2023moduli} for the details). \\

Firstly, we fix the combinatorics of the fans of the considered quantum toric stacks.

\begin{Def}
The combinatorial type of a fan $(\Delta,h : \R^n\to \R^d)$ is the poset of its cones ordered by inclusion (see as a subposet of the power set of $\{1,\ldots,n\}$ ordered by inclusion)
\end{Def}
\begin{Ex}
    The combinatorial type of the fan \[\Delta=(0,e_1,e_2,-e_1-e_2,\Cone(e_1,e_2),\Cone(e_1,-e_1-e_2),\Cone(e_2,-e_1-e_2))\]
    of $\P^2$ is 
    \[
    \{\emptyset,\{1\},\{2\},\{3\},\{1,2\},\{1,3\},\{2,3\}\}
    \]
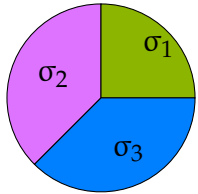
\begin{center}
    \begin{tikzpicture}[scale=.5]

 \tikzmath{\ray = 2.5;}
    \filldraw[fill opacity=0.5,fill=couleur1] (0,0) -- (0,\ray) arc (90:225:\ray cm) -- cycle;
    \filldraw[fill opacity=0.5,fill=couleur4] (0,0) -- (\ray,0) arc (0:90:\ray cm) -- cycle;
    \filldraw[fill opacity=0.5,fill=couleur3] (0,0) -- (-0.707*\ray,-0.707*\ray) arc (225:360:\ray cm) -- cycle;
    \node at (.8,.8) [above right] {$\sigma_1$};
    \node at (-2,0) [above right] {$\sigma_2$};
    \node at (0,-2) [above right] {$\sigma_3$};
\end{tikzpicture}
\end{center}
We note $S_2$ this poset (more generally, we note $S_d$ the combinatorial type of the fan of $\P^d$).
\end{Ex}

\begin{Ex}
    The combinatorial type of the fan of a complete toric surface (with $n$ generators) is 
    \[
    \{\emptyset,\{1\},\{2\},\ldots,\{n\},\{1,2\},\{2,3\},\ldots,\{n-1,n\},\{n,1\}\}
    \]
    We note $C_n$ this poset.
\end{Ex}

\begin{Def}
    Let $D$ be the combinatorial type of a fan.
A calibration $(h : \R^n \to \R^d,\Ical)$ is $D$-admissible if for all $I \in D$, $\Cone(h(e_i), i \in I)$ is strongly convex.
\end{Def}

\begin{Rem}
    By definition, a calibration $(h,\Ical)$ is $D$-admissible if there exists a quantum fan $(\Delta,h,\Ical)$ such that 
\[
\mathrm{comb}(\Delta) \simeq D
\]
Hence, the set of $D$-admissible calibrations is the set of quantum fans of combinatorial type $D$.
\end{Rem}

Now, we can define our moduli spaces of fans.

\begin{Not}
    We note $\sitename{Man}_\R$ the site of ($\Ccal^\infty$-)manifolds with smooth maps and Euclidean coverings.
\end{Not}

    \begin{Def}
    The moduli space of quantum toric stacks of dimension $d$, with $n$ generators and of combinatorial type $D$ is the stack over $\sitename{Man}_\R$
    \[
    \Mscr(d,n,D)=\{h : \R^{n} \to \R^d \mid h \text{ is } D\text{-admissible} \}/\mathrm{iso}
    \]
    (we have supposed that the combinatorial type is maximal i.e. $\Ical=[\![1,n]\!] \setminus \Delta(1)$).
    \end{Def}

Since we consider elements of some $\R^N$ and differentiable action on them, the choice of the differentiable site is natural.

\begin{Thm}[\cite{boivin2023moduli} Corollary 2.2.14 \& Proposition 2.2.15]
    Let $D$ be the combinatorial type of a complete simplicial fan. Then $\Mscr(d,n,D)$ is an orbifold given by the global quotient of a connected semi-algebraic open subset $\Omega(d,n,D)$ of $\R^{d(n-d)}$ (by fixing $d$ vectors by automorphism) through the group $\mathrm{Aut}_{\catname{Poset}}(D)$.
\end{Thm}

\begin{Thm}[\cite{boivin2023moduli} Theorem 2.4.13]
    There exists an universal family of quantum toric stacks of dimension $d$, with $n$ generators and of combinatorial type $D$ over $\Mscr(d,n,D)$ i.e. a morphism $\Xscr \to \Mscr(d,n,D)$ whose fiber are quantum toric stacks of combinatorial type $D$.
\end{Thm}

\section{Non-simplicial quantum toric stacks}
\label{NSQTV}

The paper \cite{boivin2020nonsimplicial} gives a construction of quantum toric stacks associated with an arbitrary quantum fan. This construction verify some good properties: it induces an equivalence between fans and toric stacks (see \cite[Theorem 4.2.2.2]{boivin2020nonsimplicial}), toric stacks defined in this way are global quotient stack (see \cite[Theorem 5.2]{boivin2020nonsimplicial})).\\
However, it behaves badly when we consider this construction in family: the group by which we mod out in \cite{boivin2020nonsimplicial} depends on the calibration we are considering. The aim of this chapter is to present a variant of this construction in which the group to be considered does not vary (and will be $\C^{n-d}$ as in \cref{QGIT}).
\subsection{Definitions}
\label{7-Def}

We begin with the definition of quantum toric stack associated to a arbitrary fan.
\begin{Def}
Let $(\Delta,h\colon \Z^n \to \Gamma,\Ical)$ be a quantum fan. The associated quantum toric stack $\Xscr_{\Delta,h,\Ical}$ is the quotient stack over $\sitename{A}$ \[
\Xscr_{\Delta,h,\Ical} \coloneqq [\Sscr(\Delta)/\C^{n-d}]
\]
where $\Sscr(\Delta) \subset \C^n$ is the classical toric variety described in \cref{QGIT} and the action of $\C^{n-d}$ on $\Sscr$ is defined by 
\begin{equation} \label{action_Gale_nonsimp}
    t \cdot z=E(k(t))z
\end{equation}
where $k : \R^{n-d} \to \R^,$ is a Gale transform of $h \otimes id_\R : \R^n \to \R^d$.
\end{Def}

\begin{Rem}
\begin{itemize}
    \item For simplicial fan, this definition  of quantum toric stacks coincide with the definition of quantum toric stack given in \cite{katzarkov2020quantum} (thanks to the Quantum GIT cf. \cref{QGIT})
    \item We recover the classical case by mod out the ineffectivity (in ordre to remove the gerbe) and by replacing the stacky quotient by the categorical one (see \cite[Theorem 5.1.11]{cox})
    \item The stack $\Xscr_{\Delta,h,\Ical}$ is an equivariant analytic stack (i.e. it is the stackification of an internal groupoid of $\sitename{A}$) : 
    \[
    \Xscr_{\Delta,h,\Ical}=[\pi_1,act \colon \Sscr \times \C^{n-d}\rightrightarrows \Sscr]
    \]
    where the Lie group which acts on $\Sscr \times \C^{n-d}$ is the Lie product group $((\C^*)^{n},\times) \times (\C^{n-d},+)$ and the one in $\Sscr$ is $((\C^*)^{n},\times)$.
\end{itemize}

\end{Rem}

We will recover the dense torus of $\Xscr_{\Delta,h,\Ical}$ thanks to the Quantum GIT: 

\begin{Prop}\label{isom_tore}

Let $(h : \Z^n \to \R^d,\Ical)$ be a calibration.
The stack $\left[(\C^*)^n/\C^{n-d}\right]$ and the  torus $\Tscr_{h,\Ical}$ are isomorphic (as stacks).
\end{Prop}

\begin{Cor}
The quantum torus $\Tscr_{h,\Ical}$ is dense in  $\Xscr_{\Delta,h,\Ical}$.
\end{Cor}

\subsection{Equivalence of categories}
\label{7-Equiv}

We prove the equivalence of categories between the categories of quantum fans and the one of quantum toric stacks in the case of the construction of the last subsection (like in \cref{equiv_cat}).

\begin{Lemme} \label{fonct_ctq}
The correspondence $(\Delta,h,\Ical) \mapsto \Xscr_{\Delta,h,\Ical}$ is functorial.
\end{Lemme}

\begin{proof}

Let $(L,H) \colon (\Delta,h,\Ical) \to (\Delta',h',\Ical')$ be a quantum fan morphism. By definition, the morphism $H$ induces a toric morphism $\overline{H} \colon \Sscr(\Delta) \to \Sscr(\Delta')$ (since it preserves the combinatorics). We complete the diagram \eqref{torus_morphism} with \cref{Gale_transform}:
\[\begin{tikzcd}
	0 & {\R^{n-d}} & {\R^n} & {\R^d} & 0 \\
	0 & {\R^{n'-d'}} & {\R^{n'}} & {\R^{d'}} & 0
	\arrow[from=1-1, to=1-2]
	\arrow["k", from=1-2, to=1-3]
	\arrow["h", from=1-3, to=1-4]
	\arrow[from=1-4, to=1-5]
	\arrow[from=2-1, to=2-2]
	\arrow["{k'}"', from=2-2, to=2-3]
	\arrow["{h'}"', from=2-3, to=2-4]
	\arrow[from=2-4, to=2-5]
	\arrow["L", from=1-4, to=2-4]
	\arrow["H", from=1-3, to=2-3]
	\arrow["{\widetilde{H}}", from=1-2, to=2-2]
\end{tikzcd}\]
(since for all $x \in \R^{n-d}$, $h'(H(k(x)))=L(h(k(x)))=0$ and hence there exists  (an unique, by injectivity of $k$) $z \in \R^{n'-d'}$ such that $k'(z)=H(k(x))$. We can define $\widetilde{H}(x) \coloneqq z$ and remark that the map $\widetilde{H} : \R^{n-d} \to \R^{n'-d'}$ defined this way is linear).

Note $\hcal$ the stack morphism induced by the groupoid morphism
\[\begin{tikzcd}
	{(\Sscr(\Delta) \times \C^{n-d} \rightrightarrows \Sscr(\Delta))} && {(\Sscr(\Delta') \times \C^{n'-d'} \rightrightarrows \Sscr(\Delta')) }
	\arrow["{\left(\overline{H},\widetilde{H}_\C\right)}", shift left=1, from=1-1, to=1-3]
	\arrow["{\overline{H}}"', shift right=1, from=1-1, to=1-3]
\end{tikzcd}\]
It remains to us to prove that $\hcal$ is a toric morphism. \\
In order to do this, it suffices to prove that the following diagram of groupoid morphisms commutes 
\begin{equation}\label{diagr_gpd_tor}
\begin{tikzcd}
	{[\Sscr(\Delta) \times \C^{n-d} } & {\Sscr(\Delta)]} && {[\Sscr(\Delta') \times \C^{n'-d'} } & {\Sscr(\Delta')]} \\
	{[(\C^*)^n \times \C^{n-d}} & {(\C^*)^n]} && {[(\C^*)^{n'} \times \C^{n'-d'}} & {(\C^* )^{n'} ]} \\
	{[(\C^*)^d \times \Z^{n-d}} & {(\C^*)^d]} &&{[(\C^*)^{d'} \times \Z^{n'-d'}} & {(\C^*)^{d'}]}
	\arrow[shift left=1, from=1-1, to=1-2]
	\arrow[shift right=1, from=1-1, to=1-2]
	\arrow[shift right=1, from=1-4, to=1-5]
	\arrow[shift left=1, from=1-4, to=1-5]
	\arrow["{(\overline{H},\widetilde{H}_\C)}", shift left=2, from=1-2, to=1-4]
	\arrow["{\overline{H}}"', shift right=2, from=1-2, to=1-4]
	\arrow[shift left=1, from=2-1, to=2-2]
	\arrow[shift right=1, from=2-1, to=2-2]
	\arrow[shift left=1, from=3-1, to=3-2]
	\arrow[shift right=1, from=3-1, to=3-2]
	\arrow[shift left=1, from=2-4, to=2-5]
	\arrow[shift right=1, from=2-4, to=2-5]
	\arrow[shift left=1, from=3-4, to=3-5]
	\arrow[shift right=1, from=3-4, to=3-5]
	\arrow["{(\overline{H},\widetilde{H}_\C)}", shift left=2, from=2-2, to=2-4]
	\arrow["{\overline{H}}"', shift right=2, from=2-2, to=2-4]
	\arrow[hook', from=2-2, to=1-2]
	\arrow[hook', from=2-4, to=1-4]
	\arrow[hook', from=2-1, to=1-1]
	\arrow[hook', from=2-5, to=1-5]
	\arrow["\simeq"', from=2-2, to=3-2]
	\arrow["\simeq"', from=2-4, to=3-4]
\arrow["{(\overline{L},\widetilde{H}_\C)}", shift left=2, from=3-2, to=3-4]
	\arrow["{\overline{L}}"', shift right=2, from=3-2, to=3-4]
\end{tikzcd}
\end{equation}
Hence, by stackification, we would get
\[\begin{tikzcd}
	{[\Sscr(\Delta)/\C^{n-d}]} && {[\Sscr(\Delta')/\C^{n'-d'}]} \\
	{[(\C^*)^n/\C^{n-d}]} && {[(\C^*)^{n'}/\C^{n'-d'}]} \\
	{\Tscr_{h,\Ical}} && {\Tscr_{h',\Ical'}}
	\arrow["\hcal", from=1-1, to=1-3]
	\arrow[hook', from=2-1, to=1-1]
	\arrow[hook', from=2-3, to=1-3]
	\arrow["\hcal"', from=2-1, to=2-3]
	\arrow["\simeq"', from=3-3, to=2-3]
	\arrow["\lcal"', from=3-1, to=3-3]
	\arrow["\simeq"', from=3-1, to=2-1]
\end{tikzcd}\]
where $\lcal$ is the torus morphism induced by the couple $(L,H)$.
\\
The commutativity of the diagram \eqref{diagr_gpd_tor} comes from the commutativity of the following diagram (since the last diagram is the quotient of compatible actions of $\Z^n$) 
\[\begin{tikzcd}
	{\C^n/\C^{n-d}} & {\C^{n'}/\C^{n'-d'}} \\
	{\C^d} & {\C^{d'}}
	\arrow["{[H]}", from=1-1, to=1-2]
	\arrow["{[h]}"', from=1-1, to=2-1]
	\arrow["L"', from=2-1, to=2-2]
	\arrow["{[h']}", from=1-2, to=2-2]
\end{tikzcd}\]
where $[\cdot]$ is the projection map to quotient, and the isomorphisms of \cref{isom_tore}.
\end{proof}

\begin{Prop}
Let $(\Delta,h,\Ical)$ and $(\Delta',h',\Ical')$ be quantum fans. A torus morphism $\Tscr_{h,\Ical} \to \Tscr_{h',\Ical'}$ (given by the pair $(L,H)$)  extends to a toric morphism $\Xscr_{\Delta,h,\Ical} \to \Xscr_{\Delta',h',\Ical'}$ if, and only if, $(L,H)$ is a quantum fan morphism.
\end{Prop}
\begin{proof}
The converse is given in the proof of \cref{fonct_ctq}. Suppose that the morphism can be extended. Without loss of generality, we can suppose that the two fans contain only one maximal cone i.e. $\Sscr(\Delta)=\C^I \times (\C^*)^{I^c}$ and $\Sscr(\Delta')=\C^J \times (\C^*)^J $. 
Let $z \in \C^I \times (\C^*)^{I^c}$ and $(E(z_n))_{n \in \N_{\geq 0}}$ be a sequence of $(\C^*)^n$ tending to $z$. We have to prove that the sequence $(\overline{H} E(z_n))$ tends to an element of $\C^J \times (\C^*)^{J^c}$ since it would imply that $H$ (resp. $L$) sends the generators of $\Cone(e_i, i\in I)$ (resp. $\sigma_I$) on a $\N_{\geq 0}$-linear combination of generators of $\Cone(e_j, j \in J)$ (resp of $\sigma'_J)$. \\
In order to do that, we will use the following commutative diagram:
\[
\begin{tikzcd}
 	& (\C^*)^n \arrow[dl,hookrightarrow] \arrow[rr,"{\overline{H}}"] \arrow[dd] & & (\C^*)^{n'} \arrow[dl,hookrightarrow] \arrow[dd] \\
 	\mathscr{S} \arrow[rr, crossing over,dashed] \arrow[dd] & & {\mathscr{S}'}\\
 	& \Tscr_{h,\Ical} \arrow[dl,hookrightarrow] \arrow[rr,"{\lcal}"near start] & & \Tscr_{h',\Ical'} \arrow[dl,hookrightarrow] \\
 	\Xscr_{\Delta,h,\Ical}\arrow[rr,"{\lcal}"] & & {\Xscr_{\Delta',h',\Ical'}} \arrow[from=uu, crossing over]\\
 	\end{tikzcd}
 	\]
 By hypothesis, $\lcal([z]) \in [\C^J \times (\C^*)^J/\C^{n-d}](\C)$ and $\lcal([E(z_n)])=[\overline{H}E(z_n)]$ (where $[\cdot]$ is the quotient map). We deduce that, since the action of $\C^{n-d}$ does not change the zero coordinate, the sequence $(\overline{H}(E(z_n)))$ converges in $\C^J \times (\C^*)^{J^c}$
\end{proof}

\begin{Cor} \label{Equiv_cat_nonsimp}
The functor $(\Delta,h,\Ical) \mapsto \Xscr_{\Delta,h,\Ical}$ is an equivalence between the category of quantum fans and the category of quantum toric stacks (with stack morphisms which restricts on torus morphisms as morphisms).
\end{Cor}

\subsection{Differences}
\label{7-Diff}
The last question in this section is the study of the differences between the construction of \cite{boivin2020nonsimplicial} and that of this subsection. More precisely, we will prove that the obtained stacks by these two constructions are not isomorphic.

\begin{Lemme} \label{no_isomlemme}
The groups $(\C^*,\times) \times (\Z,+)$ and $(\C,+)$ are not isomorphic.
\end{Lemme}
\begin{proof}
The torsion group of $\C^* \times \Z$ is $\mu_\infty \times \{0\}$ (where $\mu_\infty$ is the group of roots of the unity) and the torsion group of $\C$ is trivial. We deduce that these two groups are not isomorphic.
\end{proof}

\begin{Prop} \label{no_isom}
Let $h : \R^n \to \R^d$ be a linear epimorphism, $\sigma=\sigma_I$ be a cone of dimension $d$ generated by $h(e_i)$, $i \in I$ and $h_\sigma : \R^I\to \R^d$ the restriction of $h$ on $\R^I$ (which is surjective since $\sigma$ is of dimension $d$). Then the groups  $E(\ker(h_{\sigma \C})) \times \ker(h_{\sigma \mid \Z^I})$ and $\ker(h_{\sigma \C})$ are isomorphic if, and only if,  $\sigma$ is simplicial.
\end{Prop}

\begin{proof}
By the first isomorphism theorem, we have:
\[
E(\ker(h_{\sigma \C})) \simeq \ker(h_{\sigma \C})/(\ker(h_{\sigma \C}) \cap \Z^I)=\ker(h_{\sigma \C})/\ker(h_{\sigma \mid \Z^I}).
\]
The rank of $\ker(h_{\sigma \mid \Z^I})$ is equal to
\[
\mathrm{rk}(\ker(h_{\sigma \mid \Z^I}))=|I|-\mathrm{rk}(\im(h_{\sigma \mid \Z^I})) \leq |I|-d =\dim(\ker(h_{\sigma \C}))
\]
Hence, we can decompose $\ker(h_{\sigma\C})$ as follows:
\[
\ker(h_{\sigma\C})=\Vect_\C(\ker(h_{\sigma \mid \Z^I})) \oplus \C^{\mathrm{rk}(\im(h_{\sigma \mid \Z^I}))-d}
\]
We deduce that
\[
E(\ker(h_{\sigma \C})) \simeq (\C^*)^{|I|-\mathrm{rk}(\im(h_{\sigma \mid \Z^I}))} \times \C^{\mathrm{rk}(\im(h_{\sigma \mid \Z^I}))-d}
\]
We use \cref{no_isomlemme} for the conclusion (and the fact that $|I|=d$ if and only if $\sigma_I$ is simplicial).
\end{proof}

Let $h : \Z^n \to \Gamma \subset \R^d$ be a calibration and $\sigma$ be a cone of $\R^d$.
Note $\Uscr_\sigma$ the affine quantum toric stack associated to $\sigma$ in \cite{boivin2020nonsimplicial} and $\Uscr'_\sigma$ defined in this section.

\begin{Prop} \label{noGale} 

If $\sigma$ is not simplicial then the stacks $\Uscr_{\sigma}$ and $\Uscr'_\sigma$ are not isomorphic. 
\end{Prop} 

\begin{proof}

The groupoid associated to $\Uscr_\sigma$ and the groupoid associated to $\Uscr'_{\sigma}$ cannot be Morita-equivalent since their isotropy group (i.e. the stabilizer of a point for the action) cannot be isomorphic (see \cite[Theorem 4.4]{Morita}).

More precisely, the stabilizer of the action of $\Z^{N-d} \times E(\ker(h_{\sigma\C}))$ in each point of $(\C^{\widetilde{I}} \oplus 0) \times \T^J$ is $E(\ker(h_{\sigma\C})) \times \ker(h)$ and there is no point of $\C^I \times \T^{I^c}$ with a isomorphic stabilizer for the action of $\C^{n-d}$. Indeed, the stabilizer for the action of $\C^{n-d}$ are $k^{-1}(\ker(h_{\sigma\C})\cap \C^K+\ker(h))$ on a point of $\C^{\widetilde{I}} \oplus (\ker(h_{\sigma\C})  \cap \{0\}^{K}) \times \T^{I^c}$ for $\emptyset \neq K \varsubsetneq I$ and $k^{-1}(\ker(h_{\sigma\C})+\ker(h))$ on a point of $(\C^{\widetilde{I}} \oplus \{0\}) \times \T^{I^c}$. These two groups cannot be isomorphic to
$E(\ker(h_\sigma)) \times \ker(h)$ (by \cref{no_isom}).
\end{proof}

We will see this in an example:

\begin{Ex}
Consider the standard calibration $h : \Z^4 \to \R^3$ with $h(e_4)=a_1e_1-a_2e_2+a_3e_3$. 
We have seen that in $\Uscr_\sigma$ the stabilizers can be written as follows:
\begin{itemize}
    \item the stabilizer of $0_{\C^4}$ is $\Z \times E(\C(-a_1,a_2,-a_3,1))$ which is isomorphic to $\Z \times \C$  if $a_1$,$a_2$ or $a_3$ is irrational or $\Z \times \C^*$ otherwise ;
    \item the stabilizer of points of $(\C^3 \setminus \{0\} \oplus 0) \oplus 0$ is $\ker(h) \times E(\ker(h_{\sigma\C}))$ which is isomorphic to $ \C$ if $a_1$,$a_2$ or $a_3$ is irrational or $\Z \times \C^*$ otherwise;
    \item The stabilizer of $0 \oplus \C^*(-a_1,a_2,-a_3,1)$ is $\Z \times \{0\}$ ;
    \item The stabilizer of other points is $\ker(h)$ which is isomorphic to $0$ if $a_1$,$a_2$ or $a_3$ is irrational and $\Z$ otherwise.
\end{itemize}
In another hand, in $\Uscr'_\sigma$, we have: 
\begin{itemize}
    \item The stabilizer of $0_{\C^4}$ is $\C$ ; 
    \item The stabilizer of a point of $\C^I \times \{0\}$ is $\bigcap_{i \in I} a_i \Z$ (which is always a finitely generated group). In particular, if $a_1,a_2$ or $a_3$ is irrational, the stabilizer of a point of the torus is trivial and is isomorphic to $\Z$ if there are a unique non-zero coordinate. 
\end{itemize}
\end{Ex}
\section{Secondary fans}
\label{secondary_fan}
\subsection{GIT and irrational fans}

Geometric invariant theory (as developped in \cite{Mumford}) is a method for constructing quotient of a scheme by a reductive algebraic group. It is used in toric geometry to produce wall-crossings between different toric varieties.

The principal obstacle for using this formalism in quantum framework is the fact that $(\C^n,+)$ is not a reductive group: 
\begin{Lemme}
The Lie group $(\C^n,+)$ is not reductive (i.e. its unipotent radical is not trivial)
\end{Lemme}
\begin{proof}
The group $(\C^n,+)$ is unipotent thanks to the Lie group monomorphism $\iota \colon \C^n \hookrightarrow \U_{n+1}$ (where $\U_n$ is the group of upper-triangular matrices of size $n$ with a diagonal with only 1) defined by:
\[
\forall (z_1,\ldots,z_n) \in \C^n, \iota(z_1,\ldots,z_n)=\begin{pmatrix}
1 & z_1 & \ldots & z_n \\
0 & 1 & 0 & 0 \\
0 & 0 & \ddots & 0 \\
0 & 0 & & 1
\end{pmatrix}
\]
\end{proof}

A consequence of this is that the semi-stable/stable loci are not what we want to:
\begin{Ex}
Consider the action of $\C$ on $\C^2$ given by:
\[
t \cdot (z_1,z_2)=(E(t)z_1,E(\alpha t)z_2)
\]
where $\alpha \in \R_{>0} \setminus \Q_{>0}$. Consider $a \in \R$ and the character $\chi^a : t \in \C \mapsto E(at) \in \C^*$. 

Let $f$ be a $\chi^a$-equivariant holomorphic function i.e. a holomorphic map $f \colon \C^2 \to \C$ such that 
\begin{equation} \label{equiv_holom}
    \forall t\in \C, \forall (z_1,z_2) \in \C^2, f(E(t)z_1,E(\alpha t)z_2)=E(at)f(z_1,z_2)
\end{equation}

Since $f$ is of the form 
\[
\sum_{i,j \in \N_{\geq 0}} a_{ij}z_1^i z_2^j
\]
near $(0,0)$, the equality \eqref{equiv_holom} becomes, in some neighborhood $V$ of $(0,0)$
\[
\forall t\in \C, \forall (z_1,z_2) \in V, \sum_{i,j \in \N_{\geq 0}} a_{ij}E(t(i+\alpha j)) z_1^i z_2^j=E(at)\sum_{i,j \in \N_{\geq 0}} a_{ij}z_1^i z_2^j
\]
Hence, for $i,j \in \N_{\geq 0}$, $a_{ij}=0$ or $a=i+\alpha j$ i.e. $f$ is polynomial and
\[
f(z_1,z_2)=\sum_{a=i+\alpha j} a_{ij}z_1^iz_2^j
\]
Since $(1,\alpha)$ is a $\Q$-linearly independent family then the couple $(i,j)$ (when it exists) is unique:
\[
f(z_1,z_2)=a_{ij}z_1^iz_2^j
\]
Consequently, the set of $G$-semi-stables points are:
\[(\C^2)_\chi^{ss}=
\begin{cases}
\emptyset &\text{ si } a \notin \N_{\geq 0} +\N_{\geq 0} \alpha \\
(\C^*)^2 &\text{ si } a \in \N_{> 0} +\N_{> 0} \alpha \\
\C \times \C^* &\text{ si } a \in \N_{> 0} \\
\C^* \times \C &\text{ si } a \in \N_{> 0} \alpha \\
\C^2 &\text{ si } a=0 \\
\end{cases}
\]
We no more get the classical result $(\C^2)^{ss}_\chi=\C^2 \setminus \{0\}$ of the projective line. In order to have satisfying definitions of secondary fans in the quantum case, we will consider only combinatoric constructions.
\end{Ex}
\subsection{Parametrized quantum fans}

\label{10b-construction}

Let $h : \R^n \to \R^d$ be a linear epimorphism such that $h(e_i)\neq 0$ for all $i$. Let $b \in \R^n$ such that the polytope 
\[
P_b\coloneqq \{ x \in \R^d \mid \forall i \in \{1,\ldots,n\}, \left\langle x,h(e_i)\right \rangle \geq -b_i\}
\]
is of dimension $d$. The goal of this subsection is to provide a construction of a quantum fan $(\Delta_b,h,\Ical_b)$ from the datum $(h,b)$. In order to do this, we adapt the construction of the normal fan : \\
 
We begin by associate to each face of the polytope $P_b$ a cone:
\begin{Def}[\cite{cox}] \label{dualite}
Let $P \subset \R^n$ be a polytope of dimension $d$ given by:
\[
P=\bigcap_{F \preceq P,\dim(F)=d-1}\{ x \in \R^d \mid \left\langle x, u_F \right\rangle \geq -a_F\}
\]
where the intersection is indexed by the facets of $P$ and $u_F$ is a normal vector of the face $F$. \\
Then the cone associated to the face $Q$ of $P$ is 
\[
\sigma_Q \coloneqq \Cone(u_F, Q\preceq F, \dim(F)=d-1 )
\]
\end{Def}

\begin{Rem}

Everything is well-defined since for each face $F$, the couple $(u_F,a_F)$ is defined up to multiplication by a positive scalar and a cone is stable by multiplication by a positive scalar too.
\end{Rem}

As in rational case, we have the following duality statement:
\begin{Prop}[\cite{cox} Proposition 2.3.8 for the rational case] \label{dimension_dual}
Let $P \subset \R^d$ be a polytope of dimension $d$. Then, for every face $Q$ of $P$,
\[
\dim(Q)+\dim(\sigma_Q)=d
\]
Moreover, $Q \mapsto \sigma_Q$ reverse the inclusions.
\end{Prop}

We deduce that the 1-cones correspond to faces of dimension $d-1$ of the polytope or in other words, the intersection of a polytope with an hyperplane. However, the intersection of a polytope with an hyperplane is not always of maximal dimension (one can think of the intersection of the square $[0,1]^2$ with the hyperplane $x+y=2$ which is a point). 

\begin{Lemme}
The family $\Delta_b \coloneqq \{ \sigma_Q\}_{Q \preceq P_b}$ is a complete fan of $\R^d$ whose 1-cones form a subset of $\{h(e_i), 1 \leq i \leq n \}$. 
\end{Lemme}

\begin{Not}
We note $F_{i,b}$ the intersection of $P$ with the hyperplane \[H_{i,b} \coloneqq \{ x \in \R^d \mid \left\langle x,h(e_i) \right\rangle =b_i\}\]
\end{Not}
\begin{Def}
The faces $F_{i,b}$ of $P_b$ of dimension $<d-1$ is called virtual facets of the polytope $P_b$.
\end{Def}

\begin{Rem}\label{Rem_rat_GKZ}
In the rational case, one asks that the virtual facets are empty and not of dimension $<d-1$. We do not do that here because we only consider fan of maximal dimension.
\end{Rem}

Pose $\Ical_b$ the set of $i \in \{1,\ldots,n\}$ such that $F_{i,b}$ is a virtual facets of $P_\chi$, which concludes of the quantum fan $(\Delta_b,h,\Ical_b)$. 

In what follows, we consider too the polytope of $\R^n$ denoted $P_{k^\top b}$ and defined by
\[
P_{k^{\top}b} \coloneqq \{ x \in \R^n \mid \forall i, x_i \geq 0, x-b \in \ker(k^\top)\}=\{ x \in \R^n \mid \forall i, x_i \geq 0, k^\top x=k^\top b\}
\]

\begin{Lemme}
Let $\chi=k^\top b \in \R^{n-d}$. The polytope $P_{\chi}$ is the image of $P_b$ by the affine monomorphism $x \in \R^d \mapsto h^\top(x) + b \in \R^n$. This map send the facets $F_{i,b}$ on the intersections $F_{i,\chi} \coloneqq P_\chi \cap \{x_i=0\}$.
\end{Lemme}

\begin{Rem}

This last polytope is more used in the classical context of Geometric Invariant Theory since the parameter $\chi$ corresponds to a character of the Picard group (tensored by $\R$) of the associated toric variety and hence the polytope $P_\chi$ is defined in a more intrinsic than $P_b$.
\end{Rem}

In conclusion, we prove the following classification: 
\begin{Prop} \label{image_constructionEQ}
Every complete fan $(\Delta,h,\Ical)$ with a strictly convex support function $\Delta$ (i.e. a map $|\Delta| \to \R$ linear on each cone) can be realized by this construction.
\end{Prop}
\begin{proof}
We use the proof of \cite[Proposition 14.4.1]{cox} where we replace the condition of semi-projectivity by the equivalent condition of existence of a strictly convex support function (see \cite[Theorem 7.2.4]{cox}) since this is this condition which is used in the proof. Since this proof is fully combinatorial, it can be transposed to the quantum case (by deformation). It remains to verify the hypothesis of the proposition.  \\
By hypothesis, the two first points are verified and since we have removed every full-dimensional $F_{i,b}$, the third point is also verified (in the rational case, we have to removed this condition since we can have not-virtual faces with dimension $<d-1$, see \cref{Rem_rat_GKZ}).

\end{proof}
\subsection{Secondary fans}

In this subsection, we will give statements on secondary fans (and give proof if it cannot be proved in the same way as in the rational case).

\label{10b-event_second}
Let $h : \R^n \to \R^d$ be a linear epimorphism such that, for all $i$, $h(e_i)\neq 0$. Let $k$ be a Gale transform of $h$. Note $\Cone(k^\top)$ the cone
\[
\Cone(k^\top(e_1),\ldots,k^\top(e_n)) \subset \R^{n-d}
\]

\begin{Lemme}[\cite{cox} Proposition 14.3.5]
Let $\chi \in \R^{n-d}$. The following assertions are equivalent: 
\begin{itemize}
    \item $P_\chi \neq \emptyset$
    \item $\chi \in \Cone(k^\top)$
\end{itemize}
\end{Lemme}

\begin{Prop} \label{calcul_U_{adm}}
Let $\chi \in \Cone(k^\top)$. The following assertions are equivalent:
\begin{itemize}
    \item $\dim(P_\chi)=d$ ;
    \item $\chi$ is in the interior of $\Cone(k^\top)$ ; 
\end{itemize}
\end{Prop}
\begin{proof}
Suppose $\chi$ in the interior of $\Cone(k^\top(e_i),i=1..n)$. Then, thanks to the proof of $c) \Rightarrow d)$ in \cite[proposition 14.3.6]{cox}, this implies that the $F_{i,b} \coloneqq \{ x \in P_b\mid \left\langle x,h(e_i)\right\rangle=b_i\}$ are proper faces of $P_b$ for every $i$. Hence $P_b$ is of maximal dimension. \\
Conversely, suppose $P_b$ of maximal dimension and $\chi$ is in a proper face $F$ of $\Cone(k^\top(e_i),i=1..n)$.
Since $\chi$ is in $\Cone(k^\top(e_i),i=1..n)$ then 
\[
\chi=\sum_{i=1}^n \alpha_i k^\top(e_i)
\]
where $(\alpha_1,\ldots,\alpha_n) \in P_\chi$. \\
We deduce that for all $i$, $\alpha_i k^\top(e_i) \in F$ . Hence, for all $j \notin \{i \mid k^\top(e_i) \in F\} $ (which exists since $F$ is a proper face), $\alpha_j$ is zero. Consequently, for such $j$, we have $F_{j,\chi}=P_\chi$ and :
\[
F_{j,b}=P_b
\]
which is absurd since it would imply that the dimension of $P_b$ is at most the dimension of a hyperplane (since $h(e_j) \neq 0$).

\end{proof}
The following theorem gives a description of the $\chi$ for which the fan is simplicial
\begin{Thm}[\cite{cox} theorem 14.3.14] \label{generique_simple}
Let $\chi=\sum b_i \beta_i=k^\top b \in \mathrm{Int}(\Cone(\beta))$. The following assertions are equivalent:
\begin{enumerate}
    \item For every subset $I \subset \{1,\ldots,n\}$ such that $\dim\Cone(\beta)<n-d$, $\chi \notin \Cone(k^\top(e_i),i\in I)$ ("$\chi$ is generic") ; 
    \item $P_b$ is simple of dimension $d$ and the subsets \[F_{i,a}=\{ m \in P_b \mid \left\langle x,m \right\rangle=-b_i\} \] are either faces of $P_b$ or empty.
\end{enumerate}
\end{Thm}

The elements of the $\Cone(k^\top)$ encode several combinatoric types for $h$. We can decompose this cone in cones where the combinatoric is constant:

\begin{Def}
Let $(\Delta,h,\Ical)$ be a quantum fan verifying the conditions of the \cref{image_constructionEQ}. The cone GKZ associated to $(\Delta,\Ical)$, noted $\Gamma_{(\Delta,\Ical)}$ is the image of the cone
\[
\{ b \in \R^{n} \mid \exists \varphi \in \mathrm{CSF}(\Delta), \forall j \notin I, \varphi(h(e_j))=-b_j, \forall i \in I, \varphi(h(e_i)) \geq -b_i\}
\]
through $k^\top : \R^n \to \R^{n-d}$, where $\mathrm{CSF}(\Delta)$ is the set of convex support functions of $\Delta$.
\end{Def}
These cones are called GKZ for \nom{Gelfand}, \nom{Kapranov}, \nom{Zelevinsky} (since it is introduced in \cite[Chapter 7]{gelfand2009discriminants}).
\begin{Rem}
\begin{itemize}
    \item Given $b$, the associated support function is unique.
    \item The relative interior of $\Gamma_{(\Delta,\Ical)}$ is described by strictly convex $\varphi$.
\end{itemize}

\end{Rem}

\begin{Prop}[\cite{cox} Theorem 14.4.7]
The collection $\Delta_{GKZ}$ is a fan. More precisely, $\Gamma_{(\Delta,\Ical)}$ is a face of $\Gamma_{(\Delta',\Ical')}$ if, and only if, $\Delta$ refines $\Delta'$ and $\Ical' \subset \Ical$.
\end{Prop}

\begin{Def}
The maximal cones of $\Delta_{GKZ}$ are called chambers of the secondary fans and their intersection are called wall of the secondary fan. 
\end{Def}

\begin{Prop}[\cite{cox} Théorème 14.4.9]
If $\chi$ is in relative interior of $(\Gamma_{(\Delta,\Ical)})$ (i.e. in the interior of its topological closure) then the following assertions are equivalent:
\begin{itemize}
    \item $\chi$ is generic; 
    \item $\Delta$ is simplicial ; 
    \item $\Gamma_{(\Delta,\Ical)}$ is a chamber of the secondary fan.
\end{itemize}
\end{Prop}

\begin{Ex}

Consider a linear epimorphism $h:\R^5 \to \R^2$  such that 
\begin{itemize}
    \item $(h(e_3),h(e_5)) \in \Omega(\comb(\Delta_{\P^1 \times \P^1})) $,
    \item $h(e_4)\in \Omega(\comb(\Delta_{\P^2})) $
\end{itemize}
The figures \ref{fig:flip1}, \ref{fig:flip2}, \ref{fig:flip3}, \ref{fig:flip4} illustrate the differents change of combinatorics induced by the variation of the parameter $b$. In these figures, we see on the left the polytope $P_b$ and on the right the induced fan.

\begin{figure}[h!]
    \centering
    \begin{tikzpicture}
\draw[pattern={horizontal lines}, pattern color=blue,ultra thin] (0,0)--(1,0)--(1.5,1.5) -- (0,1) -- (0,0);
\draw (1.25,2.5)--(1.5,2)--(2.25,0.5) ;

  \draw[->] (6,0) -- (6,2); 
\draw[->] (6,0) -- (6+2,0); 
\draw[->] (6,0) -- (6-1.897,.63) ;
\draw[->] (6,0) -- (6.63,-1.897) ;
\draw[dashed,->] (6,0) -- (6-1.789,-.89) ;  
\end{tikzpicture}
    \caption{Fan of a Hirzebruch surface with one virtual generator}
    \label{fig:flip1}
\end{figure}
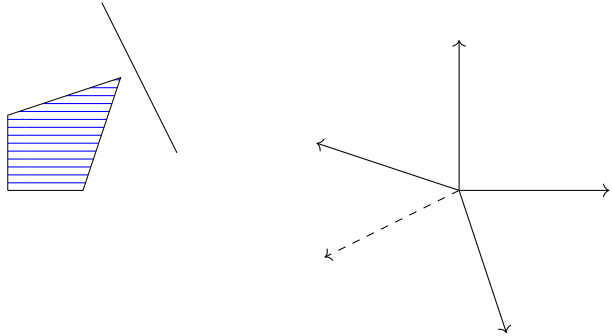

\begin{figure}[h!]
    \centering
    \begin{tikzpicture}
\draw[pattern={horizontal lines}, pattern color=blue,ultra thin] (.643,1.21) -- (0,1) -- (0,0)--(1,0)--(1.1,.3) ;
\draw (.5,1.5)--(1.25,0) ;

  \draw[->] (6,0) -- (6,2); 
\draw[->] (6,0) -- (6+2,0); 
\draw[->] (6,0) -- (6-1.897,.63) ;
\draw[->] (6,0) -- (6.63,-1.897) ;
\draw[->] (6,0) -- (6-1.789,-.89) ;  
\end{tikzpicture}
    
    \caption{Fan of the blow-up of a Hirzebruch surface in one point}
    \label{fig:flip2}
\end{figure}
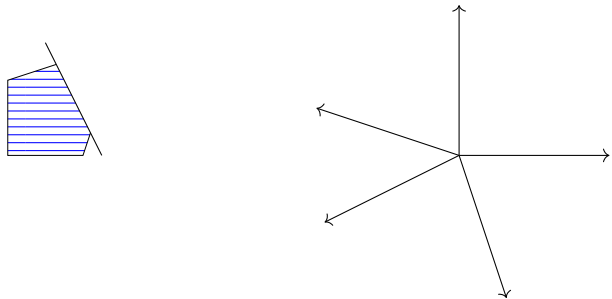

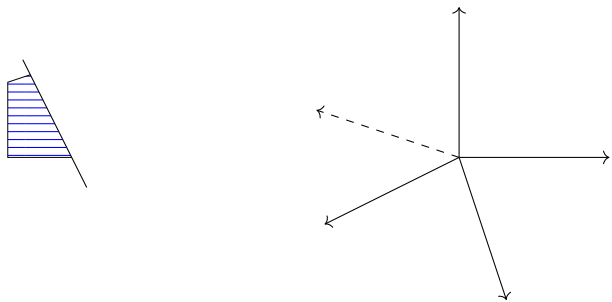
\begin{figure}[h!]
    \centering
\begin{tikzpicture}
\draw[pattern={horizontal lines}, pattern color=blue,ultra thin] (.85,0)--(0,0) -- (0,1) -- (.3,1.1);
\draw (.2,1.3)--(1.05,-.4);
  \draw[->] (6,0) -- (6,2); 
\draw[->] (6,0) -- (6+2,0); 
\draw[->,dashed] (6,0) -- (6-1.897,.63) ;
\draw[->] (6,0) -- (6.63,-1.897) ;
\draw[->] (6,0) -- (6-1.789,-.89) ;  
\end{tikzpicture}

    \caption{Fan of the blow-up of a quantum projective plane with one virtual generator}
    \label{fig:flip3}
\end{figure}

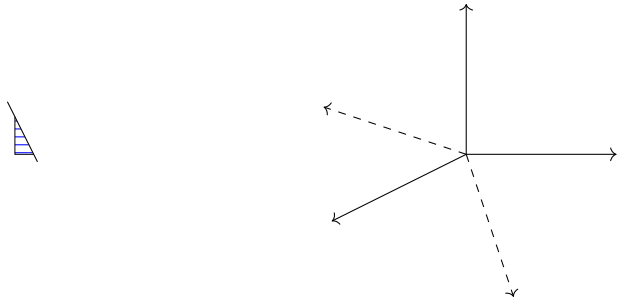
\begin{figure}[h!]
    \centering
    \begin{tikzpicture}
\draw[pattern={horizontal lines}, pattern color=blue,ultra thin] (.25,0)--(0,0) -- (0,.5);
\draw (-.1,.7)--(0.3,-.1);

  \draw[->] (6,0) -- (6,2); 
\draw[->] (6,0) -- (6+2,0); 
\draw[->,dashed] (6,0) -- (6-1.897,.63) ;
\draw[->,dashed] (6,0) -- (6.63,-1.897) ;
\draw[->] (6,0) -- (6-1.789,-.89) ;  
\end{tikzpicture}

     \caption{Fan of a quantum projective plane with two virtual generators}
    \label{fig:flip4}
\end{figure}

\end{Ex}

\subsection{Wall-crossings}
\label{10b-passagemurs}

In this subsection, we study the birational maps induced by wall-crossings i.e. between two quantum toric stacks whose combinatoric types is given by the elements in two distincts chambers separated by a wall. 

We fix a linear epimorphism $h \colon \R^n \to \R^d$ such that:
\begin{itemize}
    \item for all $i \in \{1,\ldots,n\}$, $h(e_i) \neq 0$ ; 
    \item for all $i \neq j$, $h(e_j) \notin \R h(e_i)$
\end{itemize}
Such $h$ is said geometric. 

\begin{Lemme} \label{Decomp_cone}
Let $(\Delta,h,\Ical)$ be a quantum fan given by the construction of \cref{10b-construction}. Then the associated chamber can be split in the following way:
\begin{equation}\label{Decomp_cone_formule}
\Gamma_{(\Delta,\Ical)}=\Gamma_{(\Delta,\emptyset)} \times \R_{\geq 0}^{\Ical}.
\end{equation}
\end{Lemme}

\begin{Rem} ~
\begin{itemize}
    \item Since $h$ is geometric, $\Gamma_{(\Delta,\emptyset)}$ is not empty;
    \item In the classical case, $\Gamma_{(\Delta,\emptyset)}$ is isomorphic to the cone of nef $\R$-divisors of $X_\Delta$ (included in $\Pic(X_\Delta) \otimes \R$).
\end{itemize}
\end{Rem}

\begin{proof}
Since the proof of the rational case (in \cite[Proposition 15.1.3]{cox})  is fully combinatorial (it uses support function), the equality \eqref{Decomp_cone_formule} is also true in the quantum case for chambers in the secondary fan (thanks to \cref{generique_simple} ; the case of the smaller dimension do not work since we have a stronger definition of "virtual facets" (dimension $<d-1$ instead of empty)).
\end{proof}

Hence, there are two cases for the facets of $\Gamma_{(\Delta,\Ical)}$ (i.e. the walls of secondary fans) : 
\begin{enumerate}
\item it is a product of $\Gamma_{\Delta,\emptyset} $ with $\R_{\geq 0}^{\Ical}$ ("divisorial wall") ; 
     \item it is a product of a facet of $\Gamma_{\Delta,\emptyset} $ and $\R_{\geq 0}^{\Ical \setminus \{i\}}$ where $i \in \Ical$ ("flipping wall").
\end{enumerate}

\begin{Rem}
The two adjectives associated to the wall come from the minimal model program (see, for instance, \cite{kollár_mori_1998}).
\end{Rem}

\begin{Thm} \label{qmur_renversant}
Let $(\Delta,h,\Ical)$ be a quantum fan given by the construction of \cref{10b-construction}. Let $\Gamma_{\Delta',\Ical}$ be a chamber with a wall $\Gamma_{\Delta_0,\Ical} \simeq F \times \R_{\geq 0}^\Ical$ where $F$ is a facet of $\Gamma_{\Delta,\emptyset}$. Then this is a wall between two chambers $\Gamma_{\Delta',h,\Ical}$ and $\Gamma_{\Delta,h,\Ical}$ and
\begin{itemize}
    \item $\Delta_0$ is not simplicial ;
    \item $\Delta(1)=\Delta_0(1)=\Delta'(1)$ ;
    \item $\Delta_0$ is the coarsest common refinement of $\Delta$ and $\Delta'$ ; 
    \item The exceptional loci of birational morphisms $\Xscr_{\Delta,h,\Ical} \dashrightarrow \Xscr_{\Delta_0,h,\Ical}$ et $\Xscr_{\Delta',h,\Ical}\dashrightarrow \Xscr_{\Delta_0,h,\Ical}$ are of codimension $\geq 2$ : 
    \[\begin{tikzcd}[ampersand replacement=\&]
	{\Xscr_{\Delta,h,\Ical}} \&\& {\Xscr_{\Delta',h,\Ical}} \\
	\& {\Xscr_{\Delta_0,h,\Ical}}
	\arrow[dashed,from=1-1, to=2-2]
	\arrow[dashed,from=1-3, to=2-2]
	\arrow[dashed, from=1-1, to=1-3]
\end{tikzcd}\]
\end{itemize}
\end{Thm}

\begin{proof}
The statements can be proved in the same combinatoric way as \cite[Theorem 15.3.6]{cox} (for the forth point, we can use the orbite-cone correspondence (in order to find the dimension of the  exceptional locus) in the variety $\Sscr$ and \cite[Theorem 3.3.4]{boivin2024birational}).
\end{proof}

In order to describe divisorial walls, we will need a more general notion of star subdivision:
\begin{Def}
Let $(\Delta, h \colon \Z^n \to \Gamma,\Ical)$ be a quantum fan and $v \in |\Delta| \cap \Gamma$, let $\Delta^*(v)$ the following set of cones:
\begin{itemize}
    \item $\sigma \in \Delta$ where $v \notin \sigma$ ; 
    \item $\Cone(\tau,v)$ where  $v \notin \tau \in \Delta$ and ${v} \cup \tau \subset \sigma \in  \Delta$
\end{itemize}
We call $\Delta^*(v)$ the star subdivision of $\Delta$ at $v$.
\end{Def}

Then, we have to remark that these walls corresponds, on the polytope $P_b$, to the moment when a hyperplan intersects the polytope on a smaller dimension face. Therefore, they are no more "virtual" in the classical case but they are not in the fan. In the quantum case, they remain virtual.  Hence, the toric stack in the chamber and the toric stack on the wall are the same. 
    
\begin{Thm} \label{qmur_divisoriel}
Let $(\Delta,h,\Ical)$ be a quantum fan given by the construction of \cref{10b-construction}. On the facets of $\Gamma_{\Delta,\Ical} \simeq \Gamma_{\Delta,\emptyset} \times \R^{\Ical}$ of the form $\Gamma_{\Delta,h,\emptyset} \times \R^{\Ical \setminus \{i\}}$, the obtained quantum fan is $(\Delta,h,\Ical)$ and
\begin{itemize}
    \item $\Gamma_{\Delta,\emptyset} \times \R^{\Ical \setminus \{i\}}$ is a wall between $\Gamma_{\Delta,\Ical}$ and $\Gamma_{\Delta',\Ical}$ where $\Delta'$ is the star subdivision of $\Delta$ along $\Cone(h(e_i)) \cap \Gamma$
    \item The exceptional locus of the birational map $\Xscr_{\Delta',h,\Ical} \dashrightarrow \Xscr_{\Delta,h,\Ical}$ is of codimension 1.

\end{itemize}
\end{Thm}

\begin{proof}
This is the same proof as \cite[Lemma 15.3.7]{cox} (modulo the preliminary remark and the use of \cite[Corollary 3.3.4]{boivin2024birational}).%
\end{proof}

\subsection{Flips and cobordisms}
\label{10b-cobordisme}

In this subsection, we resume the study of cobordisms (defined for quantum fans in \cite{boivin2024birational}) thanks to the secondary fan. \\
Let $P$ and $Q$ be two polytopes and $W$ be a cobordism between $P$ and $Q$ (for sake of simplicity, we suppose that the face of $W$ corresponding to $P$ is in the affine hyperplane $\{(x,t) \in \R^d \times \R \mid t=-1\}$  and the one corresponding to $Q$ is in the affine hyperplane $\{(x,t) \in \R^d \times \R \mid t=1\}$). The equations of $W$ are of the form 
\begin{equation} \label{forme1_cob}
\begin{cases}
L_1(x)+\alpha_1 t \geq -\beta_1 \\
\vdots \\
L_n(x)+\alpha_n t \geq -\beta_n \\
-1 \leq t \leq 1.
\end{cases}
\end{equation}
where the $L_i$ are linear forms $\R^d \to \R$ and $\alpha_i \in \R$. We can use also another writing of these inequations:
\begin{equation} \label{forme2_cob}
\begin{cases}
L_1(x)\geq -\beta_1-\alpha_1 t \\
\vdots \\
L_n(x)\geq -\beta_n-\alpha_n t \\
-1 \leq t \leq 1.
\end{cases}
\end{equation}
Hence, the inequations of $P$ are
\[
\begin{cases}
L_1(x)\geq -\beta_1+\alpha_1   \\
\vdots \\
L_n(x)  \geq -\beta_n +\alpha_n\\
\end{cases}
\]
and the one of $Q$ are
\[
\begin{cases}
L_1(x)\geq -\beta_1-\alpha_1 \\
\vdots \\
L_n(x)\geq -\beta_n-\alpha_n\\
\end{cases}
\]

The inequations \eqref{forme1_cob} and \eqref{forme2_cob} induce two interpretations of cobordisms in terms of secondary fans :  
\begin{itemize}
    \item the system \eqref{forme1_cob} suggests to considerate the family of calibrations $h_t : \R^n \to \R^d$ where 
    \[
    h^\top_t(x)=\left(\frac{b_i}{b_i+\alpha_i t} L_i(x)\right)_{1 \leq i \leq n}
    \]
    i.e. we fix a point in the secondary fan and we move the calibration (We make the supplementary hypothesis that $b_i+\alpha_i t > 0$ for all $t$).
    \item the system \eqref{forme2_cob} suggests to keep the calibration $h$ but to move the point of the secondary fan through the path $\chi\colon [-1,1] \to \Cone(k^\top)$ defined by:\[
    \chi(t) \coloneqq k^\top(\beta_1+\alpha_1 t,\ldots,\beta_n+\alpha_n t).
    \]
\end{itemize}

In the first case, the morphisms $k_t^\top$, given by the transpose of a Gale transform of $h_t$, varie continuously with $t$. Consequently, the chambers of the secondary fan change with $t$ (since they are described by families of $k_t^\top(e_i)$). The change of combinatorics appears when we have enough deformed $k_t^\top$ in order to ensure that $\beta=k_t^\top(\beta_1,\ldots,\beta_n)$ change of chamber.

\begin{Ex} \label{ex_blowup}
Let $h$ be the morphism $(x,y,z,w) \in \R^4 \mapsto (x-\sqrt{2}z-x,y-z-\sqrt{2}w) \in \R^2$.
Consider the cobordism (this is an example of cobordisms induced by a blow-up cf. \cite[Construction 4.2.10]{boivin2024birational}) given by the following inequalities for $(x,y,3) \in \R^2 \times [-1,1]$: 
\[\begin{cases}
x \geq 0 \\
y \geq 0 \\
\frac{4}{3+t}(-\sqrt{2}x-y) \geq -1 \\
\frac{4}{5+t}(-x-\sqrt{2}y) \geq -1
\end{cases}
\]
It is a cobordism between $W_{t=-1}=\{x \geq 0, y \geq 0, -\sqrt{2}x-y\geq -1/2, -x-\sqrt{2}y \geq -1\}$ and $W_{t=1}=\{x \geq 0, y \geq 0, -\sqrt{2}x-y\geq -1, -x-\sqrt{2}y \geq -1/2\}$.
The morphism $h_t$ ($t \in [-1,1]$) is defined by
\[
h_t(x,y,z,w)=\left(x-\frac{4\sqrt{2}}{3+t}z-\frac{4}{5+t}w,y-\frac{4}{3+t}-\frac{4\sqrt{2}}{5+t}w\right)
\]
Hence, the morphism $k_t^\top$ is defined by 
\[
k_t^\top(x,y,z,w)=\left(\frac{4}{3-t}\sqrt{2}x+\frac{4}{3-t}y+z,\frac{4}{5+t}x+\frac{4}{5+t}\sqrt{2}y+w\right)
\]
We get the following secondary fans:

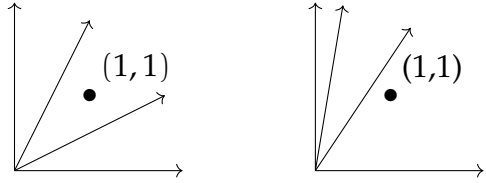
\begin{figure}[h!]
    \centering
    \begin{tikzpicture}
        \draw[->] (0,0) -- (2.236,0);
        \draw[->] (0,0) -- (1,2) ;
        \draw[->] (0,0) -- (2,1) ;
        \draw[->] (0,0) -- (0,2.236) node at (1,1){$\bullet$} node[above right] at (1,1){$(1,1)$} ; 
        \draw[->] (4,0) -- (6.236,0);
        \draw[->] (4,0) -- (4+1.1/3,2*1.1) ;
        \draw[->] (4,0) -- (4+2*1.9/3,1.9) ;
        \draw[->] (4,0) -- (4,2.236) node at (5,1){$\bullet$} node[above right] at (5,1){(1,1)} ;
    \end{tikzpicture}
    \caption{Secondary fan at $t=-1$ and at $t=1$}
    \label{fig:my_label}
\end{figure}

\end{Ex}

In the second case, the flip is an affine path in the secondary fan of the calibration $h$ between the fans $(\Delta_{\beta-\alpha},h,\Ical_{\beta-\alpha})$ and $(\Delta_{\beta+\alpha},h,\Ical_{\beta+\alpha})$. The change of combinatorics appears when we cross a wall i.e. when $t=0$. 

Conversely, an affine path of the secondary fan which crosses a wall induces a cobordism since if we have a affine path $\chi : t\in [-1,1] \mapsto k^\top(\beta+\alpha t)$ then the coefficients $\beta$ and $\alpha$ correspond to the one of \eqref{forme1_cob} and \eqref{forme2_cob}. So we have a cobordism between $(\Delta_{\chi(-1)},h,\Ical_{\chi(-1)})$ and $(\Delta_{\chi(1)},h,\Ical_{\chi(1)})$.

In conclusion, we get the following statement:
\begin{Thm}
A wall-crossing in the secondary fan given by an affine path is a cobordism. More precisely, the crossing of a divisorial fan is a cobordism of index $(1,d)$ or $(d,1)$ and the crossing of a flipping wall is a cobordism of index $(a,b)$, $a,b >1$ (see \cite[Definition 4.1.8]{boivin2024birational}). 
\end{Thm}

\begin{proof}
The first sentence is proved above. The crossing of a flipping wall has a non-simplicial "slice" (cf. \cref{qmur_renversant}) so it is a cobordism of index $(a,b)$, $a,b>1$ (cf. \cite[Warning 4.2.12]{boivin2024birational}). The crossing of a divisorial wall is of index $(1,d)$ or $(d,1)$.
\end{proof}

\subsection{Flips and projective spaces}
\label{10b-proj}

The starting point of this subsection is the previous subsection and the following statement:   
\begin{Lemme}[\cite{bosio2006} Lemma 2.3] \label{cobord_proj_polyt}
Every simple convex polytope (up to piecewise-linear homeomorphism) can be obtained from a simplex of the same dimension with a finite number of flips.
\end{Lemme}

Thanks to the previous subsection, we can ask the following natural question: can we link a quantum fan to the fan of a simplex? Or more geometrically, can we link a quantum toric stack to a quantum projective space? 

In this subsection, we will answer a bit more restrictive question: 

\begin{Quest}
    To a fixed calibration $h$, can we link a quantum fan $(\Delta,h,\Ical)$ to a fan of quantum projective space $(\Delta_{\P^d},h,\widetilde{\Ical})$? 
\end{Quest}

We pursue the study of the first question in the next section.

\begin{Thm} \label{proj_birat}
Let $h\colon \R^n \to \R^d$ be a linear epimorphism. The following statements are equivalent:
\begin{enumerate}
    \item There exists $b \in \R_{\geq 0}^n$ such that 
\[
P_b=\{ x \in \R^d \mid\forall 1 \leq i \leq n ,\left\langle  h(e_i),x\right\rangle \geq -b_i\}
\]
is a simplex;
\item There exists a subset $I\subset \{1,\ldots,n\}$ of cardinal $d+1$ such that $0 \in \Conv(h(e_i), i \in I)$.
\end{enumerate}
We note $\P(d,n)$ the set of $h$ verifying these conditions and $\P^{st}(d,n)$ the subset of $\P(d,n)$ with standard calibrations.
\end{Thm}

\begin{proof}
If there exists a $b \in \R_{\geq 0}^n$ such that $P_b$ is a simplex and there exists a set $I$ of cardinal $d+1$ such that 
\[
P_b=\{ x \in \R^d \mid \forall i \in I,\left\langle  h(e_i),x\right\rangle \geq -b_i\}
\]
We prove that $0 \in \Conv(h(e_i), i \in I)$ : \\
Let $k\in I$. By hypothesis, the restriction $h_k \coloneqq h_{\mid \R^{I \setminus \{k\}}}$ is an isomorphism. Then 
\[
h_k^\top P_b=
\{ y=h_k^\top x \in \R^{I \setminus \{k\}} \mid \forall i \in I \setminus \{k\},\left\langle  e_i,x\right\rangle \geq -b_i, \left\langle h_k^{-1}h(e_k),x\right\rangle \geq -b_k\}
\]
is also a simplex and hence $h_k^{-1}h(e_k) \in \R_{<0}^d$. We deduce that $0 \in \Conv(e_i, i \in I \setminus \{k\}, h_k^{-1}h(e_k))$. By applying $h_k$, we deduce $0 \in \Conv(h(e_i), i \in I)$. \\

Conversely, suppose that there exists a subset $I$ of cardinal $d+1$ such that $0 \in \Conv(h(e_i), i \in I)$. By the same reasoning as previously, the polytope 
\[
P\coloneqq \Conv(h(e_i), i \in I)^\circ=\{ x \in \R^d \mid \forall i \in I,\left\langle  h(e_i),x\right\rangle \geq -1\}
\]
is a simplex. We prove that there exists $b \in \R^n$ such that 
\begin{itemize}
    \item for all $i \in I$, $b_i=1$ ; 
    \item $P_b=P$.
\end{itemize}
Let $R \in \R_{>0}$ be a positive number such that the open ball $\B(0,R) \subset \R^d$ contains $P$ (which exists since $P$ is bounded thanks to the completeness of the fan). We note $R\S^{d-1}$ its boundary. 
 For all $i$, there exists $b_i$ such that the affine hyperplane
 \[
 H_{i,b_i}\coloneqq \{x \mid  \left\langle h(e_i),x \right\rangle =-b_i \}
 \]
 is tangent to $R\S^{d-1}$ and such that the upper-half plane
  \[
  H^+_{i,b_i}\coloneqq\{x \mid  \left\langle h(e_i),x \right\rangle \geq -b_i \}
 \]
 contains $B(0,R)$.\\
 Indeed, the two points where $H_{i,b_i}$ can be tangent to $R\S^{d-1}$ are 
 \[
 \pm \frac{R}{|h(e_i)|} h(e_i)
 \]
 We deduce that the possibilities for $b_i$ are 
 \[
 \pm \frac{R}{|h(e_i)|} \left\langle h(e_i),h(e_i) \right\rangle=\pm R|h(e_i)|.
 \]
 Moreover, as the upper-half plane $H^+_{i,b_i}$ contains $B(0,R)$, then 
 \[
 b_i=R|h(e_i)|.
 \]
 We found $b \in \R^n$ such that $P_b=P$ since, by construction,
 \[
 P \subset \bigcap_{i \notin I} H^+_{i,b_i}.
 \]
\end{proof}

\begin{Cor} \label{ouvert_Pdn}
The subset $\P(d,n)$ is a (non-empty) open subset of the space $\mathrm{Epi}(\R^n,\R^d)$ of epimorphisms $\R^n \to \R^d$. In other words, there exists an open subset of calibration which can be linked with a quantum projective space by the secondary fan i.e. by varying the parameter $b$.
\end{Cor}
\begin{proof}
The \cref{proj_birat} states that for all $h \in \P(d,n)$, there exists $I \subset \{1,\ldots,n\}$ such that $0 \in \Conv(h(e_i), i \in I)$. In other words, 
\[
\P(d,n)=\bigcup_{I,|I|=d+1} \{h \in \mathrm{Epi}(\R^n,\R^d) \mid 0 \in \Conv(h(e_i), i \in I) \} 
\]
If $0 \in \Conv(h(e_i), i \in I)$ then since $I$ is of cardinal $d+1$, the convex $ \Conv(h(e_i), i \in I)$ is send, by a linear transformation, on the simplex $\Conv(e_1,\ldots,e_d,v)$ où $v \in \R_{<0}^d$. We deduce that
\[
\{h \in \mathrm{Epi}(\R^n,\R^d) \mid 0 \in \Conv(h(e_i), i \in I) \} \simeq (\R^d)^{I^c} \times \GL_d(\R) \times \R_{<0}^d \subset (\R^d)^{I^c} \times (\R^d)^I
\]
which is an open subset $\mathrm{Epi}(\R^n,\R^d)$.
\end{proof}

We can do a complete description in dimension 2: 
\begin{Cor}

For all $n \geq 3$, we have two possibles cases: 
\begin{itemize}
    \item $\Omega(2,n,C_n)=\P^{st}(2,n)$ if $n \neq 4$ ; 
    \item $\Omega(2,n,C_n)\setminus \P^{st}(2,n) \simeq \R_{>0}^2$ if $n=4$
\end{itemize}
\end{Cor}
\begin{proof}
If $n=3$ then $\Omega(2,3,C_3)$ is the set of calibrations defining quantum projective planes. We deduce that 
\[
\Omega(2,3,C_3)=\P^{st}(2,3).
\]
Consider the case $n=4$. \\
Let $h \in \Omega(2,4,C_4)\setminus \P^{st}(2,4)$. Let $\alpha_3, \alpha_4 \in [0,2\pi[$ such that 
\[
\alpha_k=\mathrm{Arg}(h(e_k)),\quad k=3,4
\]
by the identification $\C \simeq \R^2$ and where $\mathrm{Arg} : \C^* \to [0,2\pi[$ is the argument function. Since $h \in \Omega(2,4,D_4)$ then 
\[
\frac{\pi}{2}<\alpha_k<2\pi
\]
Moreover, since $h \notin  \P^{st}(2,4)$ then the vectors $h(e_3)$ and $h(e_4)$ are not in $\R_{<0}^2$ and hence
\begin{equation} \label{proj12}
   \frac{\pi}{2}<\alpha_3\leq \pi, \ \frac{3\pi}{2}\leq \alpha_3<2\pi. 
\end{equation}
Since the cones generated by the $h(e_i)$ have to be strongly convex then
\begin{equation} \label{fortcvx34}
    \alpha_4-\alpha_3<\pi.
\end{equation}
Moreover $0 \notin \Conv(e_1,h(e_3),h(e_4))$ then the inequalities \eqref{proj12} give us $\alpha_3=\pi$. In the same manner, since $0 \notin \Conv(e_2,h(e_3),h(e_4))$ then $\alpha_4=\frac{3\pi}{2}$. \\
The isomorphism $\Omega(2,n,C_n)\setminus \P^{st}(2,n) \simeq \R_{>0}^2$ is $h \mapsto (|h(e_3)|,|h(e_4)|)$. \\

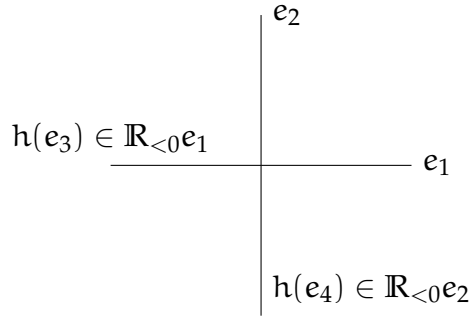
\begin{figure}[h]
\centering
\begin{tikzpicture}
    \draw (0,0) to (0,2) node[right]{$e_2$} ;
    \draw (0,0) to (2,0) node[right]{$e_1$} ;
    \draw (0,0) to (-2,0) node[above]{$h(e_3) \in \R_{<0}e_1$} ;
    \draw (0,0) to (0,-2) node[above right]{$h(e_4) \in \R_{<0}e_2$} ;
\end{tikzpicture}
    \caption{$n=4$ case}
    \label{fig:n=4case}
\end{figure}

The case $n\geq 5$ can be proved by using the case $n=4$. Suppose that there exists $h \in \Omega(2,n,C_n)\setminus \P^{st}(2,n)$. Note $k\in \{1,\ldots,n\}$ such that
\begin{itemize}
    \item $\mathrm{Arg}(h(e_i))\leq \pi$ for $i \leq k$,
    \item $\mathrm{Arg}(h(e_i))\geq \frac{3\pi}{2}$  for  $i > k$.
\end{itemize}
 Thanks to the reasoning of the case $n=4$, we have $\mathrm{Arg}(h(e_k))=\pi$ and $\mathrm{Arg} (h(e_{k+1}))=\frac{3\pi}{2}$. If $k >3$ then $0 \in \Conv(e_1,h(e_{k+1}),h(e_3))$, which is absurd. If $k=3$ then $0 \in \Conv(e_2,h(e_{k}),h(e_n))$, which is also absurd. \\ 
We can deduce that there is no element in $\Omega(2,n,C_n)\setminus \P^{st}(2,n)$.

\begin{figure}[h]
\centering
\begin{tikzpicture}
    \draw (0,0) to (0,2) node[right]{$e_2$} ;
    \draw (0,0) to (2,0) node[right]{$e_1$} ;
    \draw (0,0) to (-1.738,1.738) node[left]{$e_3$};
    \draw (0,0) to (-2,0) node[above]{$h(e_k)$} ;
    \draw (0,0) to (0,-2) node[right]{$h(e_{k+1})$} ;
    \draw[dashed] (2,0) to (-1.738,1.738) to (0,-2) to (2,0);
\end{tikzpicture}
\begin{tikzpicture}
    \draw (0,0) to (0,2) node[right]{$e_2$} ;
    \draw (0,0) to (2,0) node[right]{$e_1$} ;
    \draw (0,0) to (-2,0) node[below]{$h(e_3)$} ;
    \draw (0,0) to (0,-2) node[right]{$h(e_{4})$} ;
    \draw (0,0) to (2,-.5) node[right]{$h(e_n)$};
    \draw[dashed] (0,2) to (-2,0) to (2,-.5) to (0,2);
\end{tikzpicture}

    \caption{$n \geq 5$ case}
    \label{fig:ngeq5case}
\end{figure}
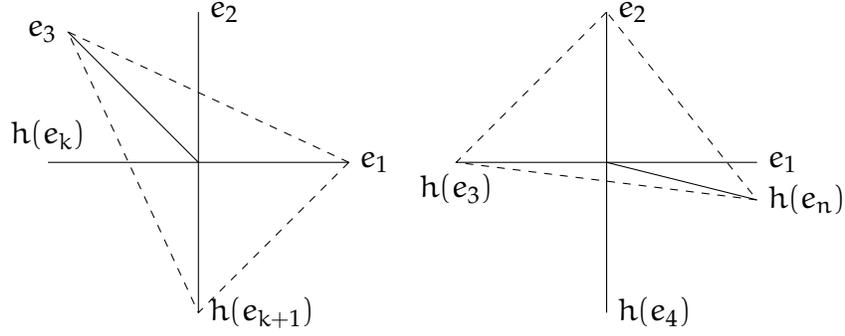
\end{proof}
\section{Augmented moduli spaces and their compactifications}
\label{augmented_ms}
Secondary fans allow us to construct paths between quantum toric stacks with different combinatorial types. In this section, we use these paths to glue together moduli spaces of different combinatorial types and their associated universal family and then, in a similar way to \cite{boivin2023moduli}, naturally compactify the resulting moduli space.
\subsection{Construction of the universal family}
\label{11-section_construction_univ}
Let $h : \R^n \to \R^d$ be a linear epimorphism. The aim of this subsection is to build a universal family of toric varieties $\Sscr$ (given by the combinatorics of the quantum fans $(\Delta_b,h,\Ical_b)$) over the secondary fan of $h$. \\
Let $b \in \R^n$ and $(\Delta_b,h,\Ical_b)$ be the associated quantum fan (cf. \cref{10b-construction} for the details). \\
The cone $\sigma=\sigma_I=\Cone(h(e_{i_1}),\ldots,h(e_{i_p}))$ is a maximal cone of $\Delta$ (i.e. of dimension $d$) if, and only if, the face corresponding to $\sigma$ in the associated  
\[P_b=\{
x \in \R^d \mid \forall i \in [\![1,n]\!], \left\langle x,h(e_i) \right\rangle \geq -b_i
\}
\] 
is the vertex given by the intersection of normal hyperplanes to the vectors $h(e_{i_j}), j \in I$. The following lemma will permit us to translate these conditions into inequations: 

\begin{Lemme}
These conditions are equivalent to the three following conditions: 
\begin{itemize}
    \item The family $(h(e_{i_k}))_{k \in [\![1,p]\!]}$ generates $\R^d$ ;  
    \item The $F_{i,b}$ are faces of dimension $d-1$ of $P_b$ (corresponding to the 1-cones of $\sigma$ in \cref{dualite}); 
    \item The equality \begin{equation}\label{dimension_intersection}\dim\{ x \in P \mid \forall i \in I, \left\langle x,h(e_i)\right\rangle = -b_i  \}=0\end{equation}
    holds.
\end{itemize}

\end{Lemme}

\begin{proof}
The first condition ensures that the cone is of the good dimension and the third ensures the existence of the vertex.
\end{proof}

The condition \eqref{dimension_intersection} is equivalent to the existence of a point $x_\sigma \in \R^d$ (cf. \cite[proposition 14.2.21]{cox}) such that
\begin{equation} \label{condition_max}
    \begin{cases}
    \forall i \in I, \left\langle x_\sigma,h(e_i)\right\rangle = -b_i \\ \forall j \notin I, \left\langle x_\sigma,h(e_j)\right\rangle \geq -b_j
    \end{cases}
\end{equation}
The equation $\eqref{condition_max}$ can be written in a more concise way: 
\[
\begin{cases}
    h_\sigma^\top (x_\sigma) = -b_I
     \\ \forall j \notin I, \left\langle x_\sigma,h(e_j)\right\rangle \geq -b_j
    \end{cases}
\]
Since $\sigma$ is of maximal dimension then $h_\sigma$ is onto and hence $h_\sigma^\top$ is one-to-one which ensures the uniqueness of the solution when it exists. 
We conclude that

\begin{Lemme} \label{S_combi}
Let  $\Sscr(P_b)$ the variety associated to the quantum fan $(\Delta_b,h,\Ical_b)$ (described in \ref{Def_S}). Then, 
\begin{equation} \label{equation_S}
 \Sscr(P_b)=\bigcup_{I} \C^I \times (\C^*)^{I^c} 
\end{equation}
where the union are indexed by the subsets $I$ of  $\{1,\ldots,n\}$ such that :
\begin{enumerate}
    \item $\dim(\sigma_I)=d$ ;
    \item There exists $x \in \R^d$ such that  $h^\top(x)=-b_I$
    and, for every $i \notin I$, $\left\langle x, h(e_i) \right\rangle \geq -b_i$ ; 
    \item $\forall i \in I, \dim(F_{i,b})=d-1$.
\end{enumerate}
\end{Lemme}

\begin{Not}
In what follows, we will note $(C)$ the three conditions of the \cref{S_combi} for the triplet $(b,I,h)$.
\end{Not}

\begin{Def}
The point $b \in  \R^n$ will be said admissible if the polytope $P_b$ is of dimension $d$. We note $\widetilde{U}^{adm}(h)$ the set of admissible $b$  
\end{Def}

When $b$ moves in $\widetilde{U}^{adm}(h) \subset \R^n$, we get a family of varieties $\Sscr$ i.e. a map $\widetilde{\Sscr}^{univ}(h) \to \widetilde{U}^{adm}$ such that for $b \in \widetilde{U}^{adm}$, $\widetilde{\Sscr}^{univ}(h)_b=\Sscr(P_b)$, by defining 
\[
\widetilde{\Sscr}^{univ}(h)=\bigcup_I \left\{(b,z) \in \R^n \times (\C^I \times (\C^*)^{I^c}) \mid (b,I,h) \text{ verifies (C)} \right\}
\]
and by taking as $\pi$ the projection $\R^n \times \C^n \to \R^n$ on the first coordinates.

Thanks to the equality \eqref{equation_S}, we have the following statement:
\begin{Lemme}
Let $b \in \R^n$. For every $m \in \ker(k^\top)$, we have:
\[
\Sscr(P_b)=\Sscr(P_{b+m})
\]
\end{Lemme}
\begin{proof}
Let $b \in \R^n$ such that $P_b \neq \emptyset$, 
and $m=-h^\top(y) \in \ker(k^\top)=\im(h^\top)$. Let $x \in P_b$. Then, for all $i$, 
\[
\left\langle x+y, h(e_i) \right\rangle =\left\langle x, h(e_i) \right\rangle+\left\langle y, h(e_i) \right\rangle=\left\langle x, h(e_i) \right\rangle+\left\langle h^\top y,e_i \right\rangle \geq -b_i-m_i
\]
This proves the inclusion $P_{b+m} \supset y+P_b$ (the translation of $P_b$ by the vector $y$). We prove in the same manner the other inclusion $y+P_b \subset P_{b+m}$. We can conclude that $y+P_b=P_{b+m}$ and hence $\Sscr(P_b)=\Sscr(P_{b+m})$.
\end{proof}

\begin{Prop}
The sets $\widetilde{\Sscr}(h)$ and $\widetilde{U}^{adm}(h)$ are invariant by the action of $\ker(h)$ (where $\ker(h)$ acts through translations). Moreover, the projection $\pi$ is $\ker(h)$-equivariant.
\end{Prop}

We deduce that the projection $\pi$ descends to the quotient by $\ker(k^\top)$ i.e. $\pi$ induces a morphism
\[
\widetilde{\Sscr}^{univ}/\ker(h) \subset \R^{n-d} \times \C^n \to U^{adm}(h) \coloneqq \widetilde{U}^{adm}(h)/\ker(k^\top) \subset \R^{n-d}
\]
since
$\R^n/\ker(k^\top)\simeq \im(k^\top)=\R^{n-d}$ by injectivity of $k$. We note $\Sscr^{univ}(h)$ the obtained family.

We have already computed $U^{adm}(h)$ in \cref{calcul_U_{adm}} :
\begin{Prop}
$U^{adm}(h)=\mathrm{Int}\Cone(k^\top(e_i),i=1..n)$
\end{Prop}

\begin{Thm} \label{univ_semi-alg}
$\Sscr^{univ}(h)$ is a semi-algebraic subset of $\R^{n-d} \times \R^{2n}$.
\end{Thm}

\begin{proof}
It suffices to prove that the different conditions of the condition $(C)$ are semi-algebraic:
\begin{itemize}
    \item the condition 1 is semi-algebraic since it is equivalent to the non-vanishing of the $d \times d$ minor of the matrix given by the $h(e_i)$, $i \in I$;
    \item the condition 2 is semi-algebraic by the quantifier elimination given by Tarski-Seidenberg theorem (see \cite[Proposition 5.2.2]{bochnak2013real}) ;
    \item the condition 3 is semi-algebraic since it is equivalent to the existence of a free family of cardinal $d-1$ (i.e. to the non-vanishing of a $(d-1) \times (d-1)$ minor) and the fact that every family of cardinal $d$ in the hyperplane $\{x\mid \left\langle x,h(e_i)\right\rangle =b_i\}$ is linearly dependent (which is also semi-algebraic by stability of semi-algebraicity by finite intersections).
\end{itemize}
\end{proof}

\begin{Ex}
We take again the notations of \cref{ex_blowup}. We have $U^{adm}(h)=\R_{\geq 0}^2$.\\
Set-theoretically, we have the equality
\begin{align*}
    \Sscr^{univ}(h)=&(\overline{\sigma_1} \cap U^{adm}) \times \C_{x,y,t}^3 \setminus \{0\} \times \C^*_z \coprod \sigma_2 \times (\C_{xz}^2 \setminus \{0\}) \times (\C_{yt}^2 \setminus \{0\}) \\
    &\coprod (\overline{\sigma_3} \cap U^{adm}) \times \C_{x,y,z}^3 \setminus \{0\} \times \C^*_t
\end{align*}
where $\overline{\sigma}$ is the (topological) closure of $\sigma$.
\end{Ex}

\subsection{Gluing of moduli spaces}
\label{11-Recollement}

In this subsection, we will do the construction of the last subsection in families indexed by $h$ in order to get a space containing a "thickening" of the open subsets $\Omega(D)$ of $D$-admissible calibrations (with D a combinatorial type).
Consider 
\[
 \widetilde{\Sscr}^{univ}=\bigcup_{I \subset \{1,\ldots,n\}} \left\{(b,z,h) \in  \R^{n-d}\times (\C^I\times (\C^*)^{I^c})\times H  \mid (b,I,h)\text{ verifies (C)} \right\} \\ 
\]
 where $H$ is the set of geometric calibrations (seen as elements of $\R^{d(n-d)}$) and ${\Sscr}^{univ}$ the image of $\widetilde{\Sscr}^{univ}$ by the map: 
\[
(b,z,h) \in \R^{n} \times \C^n \times \R^{d(n-d)} \mapsto (k^\top(b),z,h) \in \R^{n-d} \times \C^n \times \R^{d(n-d)} 
\] where $k$ is a Gale transform of $h \in \R^{d(n-d)}$.

Thanks to the previous subsection, we have
\begin{Lemme}

The projection of $\Sscr^{univ}$ on $\R^{d(n-d)}$ defines a continuous family  $(\Sscr^{univ}(h))_{h \in H}$
\end{Lemme}

 \begin{Lemme} \label{H_sa}
     The set $H$ is a semi-algebraic subset of $\R^{d(n-d)}$.
 \end{Lemme}
 \begin{proof}
     For $i,j$, the sets $H_{i,j}\coloneqq\{(v_1=e_1,\ldots,v_d=e_d,v_{d+1}\ldots,v_n) \in (\R^d)^{n} \mid v_i\notin\R v_j\}$ are semi-algebraic thanks to Tarski-Seidenberg theorem and the following equalities
     \begin{align*}H_ {ij}^c&=\{(v_1,\ldots,v_{n-d}) \in (\R^d)^{n-d} \mid \exists \lambda \in \R, v_i=\lambda v_j\}\\
     &=\pi\{(v_1,\ldots,v_{n-d},\lambda) \in (\R^d)^{n-d}\times \times \R \mid v_i=\lambda v_j\}\end{align*}
     where $\pi$ is the projection $\R^{d(n-d)} \times \R \to \R^{d(n-d)}$.
    The set $H$ is semi-algebraic as the projection of the semi-algebraic set
    \[
    \bigcap_{1 \leq i \neq j \leq n} H_{ij} \cap (\R^d \setminus \{0\})^n
    \] through the projection $(\R^d)^n \to (\R^d)^{n-d}$ forgetting the $d$ first coordinates.
 \end{proof}

In the same way as \cref{univ_semi-alg} and thanks to \cref{H_sa}, we have
\begin{Thm}
$\Sscr^{univ}$ is a semi-algebraic subset of $\R^{n-d} \times \R^{2n} \times \R^{d(n-d)}$.
\end{Thm}

It forms a family of toric varieties $\Sscr$ over its image by the projection $\R^{n-d} \times \R^{2n} \times \R^{d(n-d)} \to \R^{n-d} \times \R^{d(n-d)}$.

\begin{Lemme}
The image of the projection of $\Sscr^{univ}$ on $\R^{n-d} \times \R^{d(n-d)}$ is:
\[
U^{adm}\coloneqq \coprod_{h \in H} U^{adm}(h)=\{(\chi,h) \mid \chi \in \Cone(k^\top(e_i), i=1..n)\}
\] 
which is an open subset of $\R^{n-d} \times \R^{d(n-d)}$. 
\end{Lemme}

We have a more precise result of preservation of the combinatorics thanks to the continuity of the Gale transform
\begin{Prop} \label{ouv_pres_comb}
Let $(\chi,h_0)$ be an generic admissible couple of $\chi$. Then there exists an open neighborhood $V$ of $h_0$ in $\R^{d(n-d)}$ such that for all $h \in V$, \[\Sscr^{univ}(h)_\chi=\Sscr^{univ}(h_0)_\chi\]
\end{Prop}

\begin{proof}
The walls of the secondary fan are cones of the form $\Cone(k^\top(e_{i_1}),\ldots,k^\top(e_{i_p}))$. Since $\chi$ is generic, one can take neighborhood around the walls in order to avoid $\chi$. Consequently, by continuity of the Gale transform, there exists an open neighborhood $V$ of $h_0$ such that for every $h$ in $V$, the combinatorics of fan given by the couple do not change.
\end{proof}

\begin{Thm} \label{Thm_cnx_proj}
    For all simplicial quantum fan $(\Delta,h,\Ical)$, there exists a path in $U^{adm}$ linking it to a quantum projective space.
\end{Thm}

\begin{proof}
Let $(\Delta,h,\Ical)$ be a quantum simplicial fan. By (path-)connectedness of the $\Omega(\mathrm{comb}(\Delta))$ (see \cite[Proposition 2.2.15]{boivin2023moduli}), there exists a path between $(\Delta,h,\Ical)$ and an element $(\Delta_0,h_0,\Ical_0)$ of $\P(d,n) \cap \Omega(\mathrm{comb}(\Delta))$ (which is not empty thanks to \cref{ouvert_Pdn}). By definition of $\P(d,n)$, there exists a path in the secondary fan between $(\Delta_0,h_0,\Ical_0)$ and a quantum projective space.
\end{proof}

\begin{Rem}
Thanks to \cref{cobord_proj_polyt}, we know that every simple polytope are cobordant (up to PL homeomorphism) to a simplex. Consequently, every fan can be deformed in order to get a fan of $\P(d,n)$ which can be linked to a fan of a quantum projective space by an affine path. We can remark that we cannot use the theorem
\cite[Theorem 4.2.19]{boivin2024birational} since if the fan is not in $\P(d,n)$ then the fan obtained by deformation of the cobordism can not be in the moduli spaces of quantum projective spaces 
\end{Rem}

We now define the augmented moduli spaces which parametrize the couples $(\chi,h)$, where $\chi\in \Cone(k^\top)$ where the numbers of generators and the dimension of the ambient space is fixed, up to calibrations isomorphisms. Thanks to the equivalence of the theorems \ref{equiv_cat} and \ref{Equiv_cat_nonsimp} and to the construction of the secondary fan, it corresponds to a quantum toric stacks of fixed dimension and the fixed numbers of generators up to toric isomorphisms.

\begin{Def} \label{Def_Bdn}
The augmented moduli space $\Ascr(d,n)$ is the quotient of $U^{adm}$ by the equivalence relation  $\sim$ defined as follows: $([b],h) \sim ([b'],h')$ if $[b]=[b']$  and if there exists (for any representative) a quantum fan isomorphism between the fan induced by $h$ and $b$ and the one induced by $h'$ and $b'$ (see \cref{10b-construction}). More precisely, if we note $\catname{A}(d,n)=(X \rightrightarrows U^{adm})$ the groupoid given by this relation, $\Ascr(d,n)$ is the stackification of the pseudo-functor 
\[
T \in \sitename{Man}_\R^{op} \mapsto (X(T) \rightrightarrows U^{adm}(T)) \in \bicatname{Gpd}
\]
where $X(T)$ is the set of maps $T \to X$ which extends to a smooth map $T \to U^{adm} \times U^{adm}$ (cf. \cite[1.33]{iglesiaszemmour:hal-01288504}) and $U^{adm}(T)=\Hom_{\sitename{Man_\R}}(T,U^{adm})=\mathcal{C}^\infty(T,U^{adm})$. 
\end{Def}

\begin{Lemme} \label{presentation_edm}
The groupoid $\catname{A}(d,n)$ is equivalent to the groupoid $(R \rightrightarrows U^{adm}]$ where: 

\begin{itemize}
    \item $R \subset \GL_n(\Z)  \times U^{adm}$ is the space of $(H,b,h)$ where $H=\begin{pmatrix} P_\sigma & 0 \\ 0 & P_\tau  \end{pmatrix}$ where $(\sigma,\tau) \in \Aut(\comb(\Delta_b,h)) \times \Sfrak_{\Ical_b}$ ($(\Delta_b,\Ical_b)$ are defined in \cref{10b-construction}) ;
    \item The map source is the projection on $U^{adm}$ ; 
    \item The target map is defined as follows:
    \[
    t(H,[b],[h])=([b],L_\sigma hH^{-1})
    \]
    where $L_\sigma=(L_{\sigma(1)}  \ldots L_{\sigma(d)})^{-1}$.
\end{itemize}

\end{Lemme}
\begin{proof}
It comes from the fact that the isomorphisms of $(\Delta_b,h,\Ical_b)$ are given by permutations $(\sigma,\tau) \in \Aut(\mathrm{comb}(\Delta_b)) \times \Sfrak_{\Ical_b}$ (by \cite[Theorem 2.2.6]{boivin2023moduli}).
\end{proof}

Then we can build a universal family of associated quantum toric stacks (in a same manner as \cite[subsection 2.4]{boivin2023moduli}) by considering the groupoid $\catname{E}=(R \rightrightarrows \Sscr^{univ})$ where  
   \begin{itemize}
   \item $R \subset \GL_n(\Z) \times \C^{n-d} \times \Sscr^{univ}$ is the space of $(H,t,[b],z,h)$ where $H=\begin{pmatrix} P_\sigma & 0 \\ 0 & P_\tau  \end{pmatrix}$ where $(\sigma,\tau) \in \Aut(\comb(\Delta_b,h)) \times \Sfrak_{\Ical_b}$ where $(\Delta_b,\Ical_b)$ are defined in \cref{10b-construction} ;
    \item The source map is the projection on $\Sscr^{univ}$ ; 
    \item The target map is defined as follows
    \[
    t(H,t,[b],z,[h])=([b],E(Hk\widehat{H}^{-1}t)z,L_\sigma hH^{-1})
    \]
where $\widehat{H}$ is the linear automorphism of $\C^{n-d}$ which makes the following diagram commutes
\[\begin{tikzcd}[ampersand replacement=\&]
	{\C^{n-d}} \& {\C^{n-d}} \\
	{\C^n} \& {\C^n}
	\arrow["{\widehat{H}}", from=1-1, to=1-2]
	\arrow["k"', from=1-1, to=2-1]
	\arrow["{k'}", from=1-2, to=2-2]
	\arrow["H"', from=2-1, to=2-2]
\end{tikzcd}\]
   \end{itemize}
   
We will note $\Escr(d,n)$ the stack obtained by stackification of $\catname{E}$. The projection induces a morphism $\Escr(d,n) \to \Ascr(d,n)$ whose fibers are quantum toric stacks. It is the universal family over the augmented moduli space $\Ascr(d,n)$. 
Thanks to this presentation, we deduce the following. 
\begin{Lemme} \label{top_stacks}
The stacks $\Ascr(d,n)$ and $\Escr(d,n)$ are topological stacks in the sense of \cite{Noohi2005FoundationsOT} (if we forget the differentiable structure).%
\end{Lemme}

\begin{Thm} \label{equiv_homot}
Let $n \geq d$ two integers. Let $D$ be the combinatorial type of a simplicial fan verifying the conditions of \cref{image_constructionEQ} (and $\Ical=D(1)^c$ be the associated set of virtual generators). Let $(\Delta,h,\Ical)$ be a quantum fan of combinatorial type $D$. Note $\widehat{\Omega}(D)$ the set of couples $(h,[b])$ of $\Omega(D) \times \R^{n-d} \cap U^{adm}$ which is in the same chamber of the secondary fan as$(\Delta,h,\Ical)$ (i.e. which has the combinatorics as $(\Delta,h,\Ical)$). Then $\widehat{\Omega}(D)$ is an open subset of $U^{adm}$ and the projection \[\widehat{\Omega}(D) \to \Omega(D)\] is a $\R^{n-d}$-fibration and hence a homotopy equivalence.
\end{Thm}

\begin{Rem}
The combinatorial type $D$ can be not simplicial but then $\widehat{\Omega}(D)$ is not open in $U^{adm}$
\end{Rem}

\begin{proof}
The set $\widehat{\Omega}(D)$ is an open subset of $\Omega(D) \times \R^{n-d} \cap U^{adm}$ and thus of $U^{adm}$ by \cref{ouv_pres_comb}. Each fiber of the projection $\widehat{\Omega}(D) \to \Omega(D)$ is an open convex cone and hence is homeomorphic to  $\R^{n-d}$. By continuity of the Gale transform, it is a fibration of fiber $\R^{n-d}$ which is contractile. The zero section is the desired quasi-inverse.
\end{proof}

In the same way, $\Sscr^{univ}_{\widehat{\Omega}(D)} \to \Omega(D) \times \Sscr(D)$ is a $\R^{n-d}$-fibration too (where $\Sscr^{univ}_{\widehat{\Omega}(D)}$ is the restriction de $\Sscr^{univ}$ sur $\widehat{\Omega}(D)$). 
\\
In order to get a family of quantum toric stacks, we can make act $\C^{n-d}$ on $\Sscr^{univ}$ by: 
\begin{equation} \label{action_Cn-d_Omegachap}
   t \cdot ([b],z,h)=([b],E(k(t))z,h) 
\end{equation}
The quotient stack $[\Sscr^{univ}/\C^{n-d}]$ will be denoted $\Xscr^{univ}$ and

\begin{Cor}
The projection $\widehat{\Xscr}(D) \coloneqq \left[\Sscr_{\widehat{\Omega}(D)}/\C^{n-d}\right] \subset \Xscr^{univ} \to \Xscr(D)$ is a $\R^{n-d}$-fibration where $\Xscr(D)$ is the total space of a family of quantum toric stacks of combinatorial type $D$ over $\Omega(D)$.
\end{Cor}

Then, we can define an action of $G \coloneqq \mathrm{Aut}(D) \times \Sfrak_\Ical \subset \Sfrak_n$ on $\widehat{\Xscr}(D)$ by descent to quotient of the action on $\widehat{\Omega}(D)$ by
\[
\sigma \cdot ([b],z,[k])=\left([b_{\sigma}],z_\sigma,\sigma \cdot h \right)
\]
where $b_{\sigma}$ is the point $(b_{\sigma(1)},\ldots,b_{\sigma(n-d)})$. \\
In the same way, we can define an action on  $\widehat{\Omega}(D)$ by: 
\begin{equation} \label{action_G_chap}
    \sigma \cdot ([b],h)=\left([b_{\sigma}],\sigma \cdot h \right)
\end{equation}

\begin{Rem}
This action is natural because if $(L,H)$ is an isomorphism of quantum fans between $(\Delta,h,\Ical)$ and  $(\Delta',h',\Ical')$ and \[
P\coloneqq \{x \in (\R^d)^\vee \mid \left\langle x,h'(e_i)\right\rangle \geq -b_{\sigma(i)} \}
\]is the polytope associated to $(\Delta',h',\Ical')$ in the dual space $(\R^d)^\vee$ then
\begin{align*}
    L^\top P&=\left\{ L^\top x \mid \left\langle x,h'(e_i)\right\rangle \geq -b_{\sigma(i)}\right\} \\
    &=\left\{ y\mid \left\langle (L^\top)^{-1}(y),h'(e_i)\right\rangle \geq -b_{\sigma(i)}\right\} \\
    &=\left\{ y\mid \left\langle y,L^{-1}h'(e_i)\right\rangle \geq -b_{\sigma(i)}\right\} \\
     &=\left\{ y\mid \left\langle y,hH^{-1}(e_i)\right\rangle \geq -b_{\sigma(i)}\right\} \\
     &=\left\{ y\mid \left\langle y,h(e_{\sigma(i)})\right\rangle \geq -b_{\sigma(i)}\right\} \\
      &=\left\{ y\mid \left\langle y,h(e_i)\right\rangle \geq -b_{i}\right\} \\
\end{align*}
\end{Rem}
These actions make the projections $\widehat{\Xscr}(D) \to \widehat{\Omega}(D)$, $\widehat{\Xscr}(D) \to \Xscr(D)$ and $\widehat{\Omega}(D) \to \Omega(D)$  $G$-equivariant maps.
They descend to quotient : 
\[\pcal_1 : \widehat{\Omega}(D)/G \to [\Omega(D)/G]=\Mscr(d,n,D), \]
\[\pcal_2 : \widehat{\Xscr}(D)/G \to \Xscr(D)/G\]
and
\[\pcal_3 : \widehat{\Xscr}(D)/G \to [\widehat{\Omega}(D)/G]\]

\begin{Cor} \label{Cor1}
The morphism $\pcal_1 : \widehat{\Omega}(D)/G \to \Mscr(d,n,D)$ is a $\R^{n-d}$-fibration.
\end{Cor}

\begin{proof}
By \cref{equiv_homot}, the composition $\widehat{\Xscr}(D) \to \Omega(D) \to [\Omega(D)/G]$ is a $G \times \R^{n-d}$-fibration. We deduce that the morphism obtained by descent to quotient $\widehat{\Omega}(D)/G \to \Omega(D)/G$ is a $\R^{n-d}$-fibration. 
\end{proof}

We have the same statement for the morphism $\pcal_2$ :
\begin{Cor} \label{equiv_homot_p2}
The morphism $\pcal_2 : \widehat{\Xscr}(D)/G \to \Xscr(D)/G$ is a $\R^{n-d}$-fibration.
\end{Cor}

In conclusion, we have

\begin{Thm}
    Let $d \leq n$ two integers. For every quantum fan (verifying the conditions of \cref{image_constructionEQ}) in $\R^d$ with $n$ generators and the combinatorial type $D$, we have the following statements:
    \begin{itemize}
        \item there exists a substack $\Ascr(D)$ (which is open if $D$ is simplicial) of $\Ascr(d,n)$ and a $\R^{n-d}$-fibration $\Ascr(D) \to \Mscr(d,n,D)$ ;
        \item there exists a substack $\Escr(D)$ (which is open if $D$ is simplicial) of $\Escr(d,n)$ and a $\R^{n-d}$-fibration $\Escr(D) \to \Xscr(D)$ ;
        \item we have a projection $\Escr(D) \to \Ascr(D)$ whose fibers are quantum toric stacks of combinatorial type $D$.
    \end{itemize}
    
\end{Thm}

\begin{proof}
We have 
\begin{itemize}
    \item $\Ascr(D)\coloneqq [\widehat{\Omega}(D)/G]$ and the projection is given by the projection of \cref{Cor1}; 
    \item $\Escr(D)\coloneqq 
    \widehat{\Xscr}(D)/G$ and the fibration is the projection of \cref{equiv_homot_p2} ; 
    \item The projection $\Escr(D) \to \Ascr(D)$ is the morphism $\pcal_3$. Its fibers are quantum toric stacks of combinatoric type $D$ by contruction (this is the descent to quotient of the restriction of the projection $\Xscr^{univ} \to U^{adm}$ to quantum fans of combinatorial type $D$).
\end{itemize}
\end{proof}
To sum up, we have the following commutative diagram
\begin{equation} \label{resume_Thm_edm_aug}
\begin{tikzcd}[ampersand replacement=\&]
	{\Xscr(d,n,D)} \&\& {\Escr(D)} \& {\Escr(d,n)} \\
	{\Mscr(d,n,D)} \&\& {\Ascr(D)} \& {\Ascr(d,n)}
	\arrow[from=1-1, to=2-1]
	\arrow["{\R^{n-d}\text{-fibration}}", from=1-3, to=1-1]
	\arrow["{\R^{n-d}\text{-fibration}}"', from=2-3, to=2-1]
	\arrow[from=1-3, to=2-3]
	\arrow[hook, from=1-3, to=1-4]
	\arrow[hook, from=2-3, to=2-4]
	\arrow[from=1-4, to=2-4]
\end{tikzcd}
\end{equation}
where the fibers of the vertical morphisms are quantum toric stacks.

\subsection{Compactification of augmented moduli spaces}
\label{11-Compactification}
This subsection will be devoted to the compactification of the augmented moduli spaces. The first thing to remark is that for every cone of the secondary fan, we can restrict us to its intersection with the closed unit ball since the combinatorics do not change in this cone.
 We will consider the closed subset $U^{adm} \cap \overline{\B}(0,1) \times \R^{d(n-d)}$ de $U^{adm}$ (where $\overline{\B}(0,1)$ is the closed unit ball of $\R^{n-d}$ for Euclidean norm). We have:

\begin{Lemme}
The inclusion $U^{adm} \cap \overline{\B}(0,1) \times \R^{d(n-d)} \hookrightarrow U^{adm}$ is a homotopy equivalence.
\end{Lemme}

\begin{proof}
 It is the consequence of that fact that for any strongly convex cone $\sigma$, $\sigma$ and $\sigma \cap \overline{\B}(0,1)$ are homotopic to $\{0\}$.
\end{proof}
Moreover, this restriction is compactible with the group actions:
\begin{Lemme}
The closed subset $\widehat{\Omega}(D) \cap \overline{\B}(0,1) \times \R^{d(n-d)}$ is invariant by the $G$-action induced by \eqref{action_G_chap}.
\end{Lemme}

\begin{proof}
Since the Euclidean norm is invariant by the action of the symmetric group, the unit ball is invariant by the action of this group by permutation of coordinates.
\end{proof}
In the same manner as \cite[section 3]{boivin2023moduli}, we can embed  $U^{adm} \cap \overline{\B}(0,1) \times \R^{d(n-d)} $ in the compact space $\overline{\B}(0,1) \times \Gr(n-d,\R^n)$. We will note $U$ the image of this embedding in $\overline{\B}(0,1) \times \Gr(n-d,\R^n)$.\\
We can transfer the action \eqref{action_Cn-d_Omegachap} to $U$.
In the same manner, the family of quantum toric stacks $U^{adm} \cap \overline{\B}(0,1) \times \R^{d(n-d)}$ can be push-forward into a family $\Fscr \to U$.

Note $K$ the closure of $U$ in $\overline{\B}(0,1) \times \Gr(n-d,\R^n)$. In a analogous way as \cite[subsection 3.2]{boivin2023moduli}, we get a family of quantum toric stacks $\overline{\Fscr} \to K$. \\

In the same manner as \cref{equiv_homot}, we have the following statement
\begin{Lemme} \label{widehat}
For every combinatorial type $D$ verifying conditions of the \cref{image_constructionEQ}, the projection $\overline{\widehat{\Omega}(D) \cap \overline{\B}(0,1) \times \Gr(n-d,\R^n)} \subset K \to \overline{\Omega}(D)$ is a
$\R^{n-d}$-fibration.
\end{Lemme}

We can define on $K$ the equivalence relation  $\sim$ defined by
$(h,[b]) \sim (k,[b'])$ si $[b]=[b']$  and if there exists (for any representative) a fan isomorphism between the fan induced by $[h]$ and $b$ between the fan indced by $[k]$ et $b'$ (see \cref{10b-construction}). This equivalence relation defines a groupoid $\catname{K}(d,n)$. We build in the same manner as \cref{Def_Bdn} a associated stack $\Kscr(d,n)$  which is the "compactification "\footnote{is a not a true compactification since we restrict to a closed unit ball} of $\Ascr(d,n)$.

As in lemmas \ref{presentation_edm} and \ref{top_stacks}, we have the following results:

\begin{Lemme}
The groupoid $\catname{K}(d,n)$ is equivalent to the groupoid $(R \rightrightarrows K)$ where: 
\begin{itemize}
    \item $R \subset \GL_n(\Z) \times K$ is the space $(H,[b],h)$ where $H=\begin{pmatrix} P_\sigma & 0 \\ 0 & P_\tau  \end{pmatrix}$ where $(\sigma,\tau) \in \Aut(\comb(\Delta_b,h)) \times \Sfrak_{\Ical_b}$ where $(\Delta_b,\Ical_b)$ are defined in \cref{10b-construction} ;
    \item The source map is the projection on $K$ ; 
    \item The target map is defined as follows
    \[
    t(H,[b],[h])=([b],[hH^{-1}])
    \]
\end{itemize}

\end{Lemme}

\begin{Lemme}
The stack $\Kscr(d,n)$ is a topological stack (if
we forget the differentiable structure).
\end{Lemme}

In the same manner as \cref{equiv_homot_p2}, we get
\begin{Lemme}
The projection of \cref{widehat} descend as a $\R^{n-d}$-fibration
\[
\Kscr(D)\coloneqq \left(\ \overline{\widehat{\Omega}(D) \cap \overline{\B}(0,1) \times \Gr(n-d,\R^n)} / \sim \right) \subset \Kscr(d,n) \to \overline{\Mscr(d,n,D)}.
\]
\end{Lemme}

\begin{proof}
It suffices to prove that fan isomorphisms of elements in $\widehat{\Omega}(D)$ are in $\Aut(D) \times \Sfrak_{D(1)^c}$.
\end{proof}
By the same method as \cite[subsection 3.2]{boivin2023moduli} or in the definition of $\Escr(d,n)$ in the previous subsection, we can push-forward the family $\Fcal \to K$ to $\Kscr(d,n)$. We will note $\overline{\Xscr(d,n)}$ the obtained total space which is the universal family over the compactification $\Kscr(d,n)$. In conclusion, we get an answer to the "Toric Big Moduli Conjecture": 

\begin{Thm}\label{Compactification_edm_recolle}
    Let $n \geq d$ two integers. There exists a compact stack (over the site $\sitename{Man_\R}$)  $\Kscr(d,n)$ and a stack morphism $\overline{\Escr(d,n)} \to \Kscr(d,n)$ such that:
    \begin{itemize}
        \item For every combinatorial type $D$ with a complete fan of $\R^d$ with $n$ generators, there exists a closed substack $\Kscr(D)$ of $\Kscr(d,n)$ and a $\R^{n-d}$-fibration $\Kscr(D) \to \overline{\Mscr(d,n,D)}$ ;
        \item This fibration can be extended to a fibration between the family over $\Kscr(D)$ and the one over $\overline{\Mscr(d,n,D)}$. 
    \end{itemize}
    \end{Thm}
We can sum up this statement by the following diagram which is the "compactification" version of the diagram \eqref{resume_Thm_edm_aug} : 

\[\begin{tikzcd}[ampersand replacement=\&]
	{\overline{\Xscr}(d,n,D)} \&\& {\overline{\Escr}(D)} \& {\overline{\Escr}(d,n)} \\
	{\overline{\Mscr}(d,n,D)} \&\& {\Kscr(D)} \& {\Kscr(d,n)}
	\arrow[from=1-1, to=2-1]
	\arrow["{\R^{n-d}\text{-fibration}}", from=1-3, to=1-1]
	\arrow["{\R^{n-d}\text{-fibration}}"', from=2-3, to=2-1]
	\arrow[from=1-3, to=2-3]
	\arrow[hook, from=1-3, to=1-4]
	\arrow[hook, from=2-3, to=2-4]
	\arrow[from=1-4, to=2-4]
\end{tikzcd}\]

In other words, the (compact) stack $\Kscr(d,n)$ contain a thickening (given by a $\R^{n-d}$-fibration) of each compactification $\Mscr(d,n,D)$ of moduli spaces of quantum toric stacks, gluing thanks to the secondary fan.


\bibliographystyle{alpha} 
\bibliography{Biblio}

\end{document}